\documentclass[aop,MSNbibl,seceqn,dvips]{arximspdf}
\usepackage{accents}

\doi{10.1214/12-AOP788} 
\volume{42}
\issue{1}
\pubyear{2014}
\firstpage{204}
\lastpage{236}

\makeatletter

\renewcommand{\underline}{\underbar}
\renewcommand{\overline}{\bar}
\renewcommand{\underbar}{\underaccent{\bar}}

\newcommand{\cal}{\mathcal}

\def\dbD{\mathbb{D}}
\def\dbE{\mathbb{E}}
\def\dbF{\mathbb{F}}
\def\dbP{\mathbb{P}}
\def\dbR{\mathbb{R}}

\def\a{\alpha}
\def\b{\beta}
\def\d{\delta}
\def\e{\varepsilon}

\def\t{\tau}
\def\f{\varphi}
\def\th{\theta}
\def\G{\Gamma}
\def\L{\Lambda}
\def\O{\Omega}

\def\cA{{\cal A}}
\def\cD{{\cal D}}
\def\cE{{\cal E}}
\def\cF{{\cal F}}
\def\cL{{\cal L}}
\def\cT{{\cal T}}
\def\cU{{\cal U}}
\def\cY{{\cal Y}}
\def\cZ{{\cal Z}}

\def\one{{\mathbf1}}

\newtheorem{theorem}{Theorem}[section]
\newtheorem{lem}[theorem]{Lemma}
\newtheorem{prop}[theorem]{Proposition}

\newproclaim{rem}[theorem]{Remark}
\newtheorem{eg}[theorem]{Example}
\newproclaim{defn}[theorem]{Definition}
\newproclaim{assum}[theorem]{Assumption}

\makeatother

\begin{document}
\begin{frontmatter}

\title{On viscosity solutions of path dependent PDEs}
\runtitle{Viscosity solutions of PPDEs}

\begin{aug}
\author[A]{\fnms{Ibrahim} \snm{Ekren}\ead[label=e1]{ekren@usc.edu}\thanksref{tt1}},
\author[A]{\fnms{Christian} \snm{Keller}\ead[label=e2]{kellerch@usc.edu}},
\author[B]{\fnms{Nizar} \snm{Touzi}\thanksref{t1,tt2}\ead[label=e3]{nizar.touzi@polytechnique.edu}}\\
\and
\author[A]{\fnms{Jianfeng} \snm{Zhang}\corref{}\thanksref{t2}\ead[label=e4]{jianfenz@usc.edu}}
\runauthor{Ekren, Keller, Touzi and Zhang}
\affiliation{University of Southern California\thanksmark{tt1} and
Ecole Polytechnique\thanksmark{tt2}}
\address[A]{I. Ekren\\
C. Keller\\
J. Zhang\\
Department of Mathematics\\
University of Southern California\\
3620 S. Vermont Ave, KAP 108\\
Los Angeles, California 90089\\
USA\\
\printead{e1}\\
\hphantom{E-mail: }\printead*{e2}\\
\hphantom{E-mail: }\printead*{e4}} 
\address[B]{N. Touzi\\
CMAP, Ecole Polytechnique\\
91128 Palaiseau Cedex, Paris\\
France\\
\printead{e3}}
\end{aug}

\thankstext[1]{t1}{Supported by the Chair \textit{Financial Risks} of
the \textit{Risk Foundation} sponsored by Soci\'et\'e
G\'en\'erale, the Chair \textit{Derivatives of the Future}
sponsored by the F\'ed\'eration Bancaire Fran\c{c}aise, and
the Chair \textit{Finance and Sustainable Development} sponsored by
EDF and Calyon.}

\thankstext[2]{t2}{Supported in part by NSF Grant DMS-10-08873.}

\received{\smonth{11} \syear{2011}}
\revised{\smonth{6} \syear{2012}}

%
\begin{abstract}
In this paper we propose a notion of viscosity solutions for path
dependent semi-linear parabolic PDEs. This can also be viewed as
viscosity solutions of non-Markovian backward SDEs, and thus extends
the well-known nonlinear Feynman--Kac formula to non-Markovian case.
We shall prove the existence, uniqueness, stability and comparison
principle for the viscosity solutions. The key ingredient of our
approach is a functional It\^{o} calculus recently introduced by Dupire
[Functional It\^o calculus (2009) Preprint].
\end{abstract}

%

\pdfkeywords{35D40, 35K10, 60H10, 60H30,
Path dependent PDEs, backward SDEs,
functional Ito formula,
viscosity solutions,
comparison principle}

\end{frontmatter}

\section{Introduction}
\label{sect-Introduction}

It is well known that a Markovian type backward SDE (BSDE, for short)
is associated with a semi-linear parabolic PDE via the so called
nonlinear Feynman--Kac formula; see Pardoux and Peng \cite{PP2}. Such
relation was extended to forward--backward SDEs (FBSDE, for short) and
quasi-linear PDEs; see, for example, Ma, Protter and Yong \cite{MPY},
Pardoux and Tang \cite{PT} and Ma, Zhang and Zheng \cite{MZZ}, and
second order BSDEs (2BSDEs, for short) and~fully nonlinear PDEs; see,
for example, Cheridito et al. \cite{CSTV} and Soner, Touzi and
Zhang~\cite{STZ-2BSDE}. The notable notion $G$-expectation, proposed by
Peng \cite{Peng-G}, was also motivated from connection with fully
nonlinear PDEs.

In non-Markovian case, the BSDEs (and FBSDEs, 2BSDEs) become path
dependent. Due to its connection with PDE in Markovian case, it has
long been discussed that general BSDEs can also be viewed as a PDE. In
particular, in his ICM 2010 lecture, Peng \cite{Peng-ICM} proposed the
question whether or not a non-Markovian BSDE can be viewed as a
path-dependent PDE (PPDE, for short).

The recent work Dupire \cite{Dupire}, which was further extended by
Cont and Fournie \cite{CF}, provides a convenient framework for this
problem. Dupire introduces the notion of horizontal derivative (that we
will refer to as time derivative) and vertical derivative (that we will
refer to as space derivative) for nonanticipative stochastic
processes. One remarkable result is the functional It\^{o} formula
under his definition. As a direct consequence, if $u(t,\omega)$ is a
martingale under the Wiener measure with enough regularity (under their
sense), then its drift part from the It\^{o} formula vanishes, and thus
it is a classical solution to the following path-dependent heat equation:
%
\begin{equation}
\label{heat} \partial_t u(t,\omega) + \tfrac12 \,\partial^2_{\omega\omega}
u(t,\omega) = 0.
\end{equation}

It is then very natural to view BSDEs as semi-linear PPDEs, and 2BSDEs
and $G$-martingales as fully nonlinear PPDEs. However, we shall
emphasize that PPDEs can rarely have classical solutions, even for heat
equations. We refer to Peng and Wang \cite{PW} for some sufficient
conditions under which a semi-linear PPDE admits a classical solution.

The present work was largely stimulated by Peng's recent paper \cite
{Peng-viscosity}, which appeared while our investigation of the problem
was in an early stage. Peng proposes a notion of viscosity solutions
for PPDEs on c\`{a}dl\`{a}g paths using compactness arguments. However, the
horizontal derivative (or time derivative) in \cite{Peng-viscosity} is
defined differently from that in Dupire \cite{Dupire} which leads to a
different context than ours. Moreover, Peng \cite{Peng-viscosity}
derives a uniqueness result for PPDEs on c\`{a}dl\`{a}g paths. Given the
nonuniqueness of extension of a function to the c\`{a}dl\`{a}g paths, this
does not imply any uniqueness statement in the space of continuous
paths. For this reason, our approach uses an alternative definition
than that of Peng \cite{Peng-viscosity}.

The main objective of this paper is to propose a notion of viscosity
solutions of PPDEs on the space of continuous paths. To focus on the
main idea, we focus on the semi-linear case and leave the fully
nonlinear case for future study. We shall prove existence, uniqueness,
stability, and comparison principle for viscosity solutions.

The theory of viscosity solutions for standard PDEs has been well
developed. We refer to the classical references Crandall, Ishii and
Lions \cite{CIL} and Fleming and Soner \cite{FS}. As is well
understood, in path-dependent case the main challenge comes from the
fact that the space variable is infinite dimensional and thus lacks
compactness. Our context does not fall into the framework of Lions
\cite{Lions-Hilbert-I,Lions-Hilbert-II,Lions-Hilbert-III} where the
notion of viscosity solutions is extended to Hilbert spaces by using a
limiting argument based on the existence of a countable basis.
Consequently, the standard techniques for the comparison principle,
which rely heavily on the compactness arguments, fail in our context.
We shall remark though, for first order PPDEs, by using its special
structure Lukoyanov \cite{Lukoyanov} studied viscosity solutions by
adapting elegantly the compactness arguments.

To overcome this difficulty, we provide a new approach by decomposing
the proof of the comparison principle into two steps. We first prove a
partial comparison principle, that is, a classical sub-solution (resp.,
viscosity sub-solution) is always less than or equal to a viscosity
super-solution (resp., classical super-solution). The main idea is to
use the classical one to construct a test function for the viscosity
one and then obtain a contradiction.

Our second step is a variation of the Perron's method. Let $\underline
u$ and $\overline u$ denote the supremum of classical sub-solutions and
the infimum of classical super-solutions, respectively, with the same
terminal condition. In standard Perron's approach (see, e.g., Ishii
\cite{Ishii} and an interesting recent development by Bayraktar and
Sirbu \cite{BS}), one shows that
%
\begin{equation}
\label{u=u} \underline u=\overline u
\end{equation}
by assuming the comparison principle for viscosity solutions, which
further implies the existence of viscosity solutions. We shall instead
prove (\ref{u=u}) directly, which, together with our partial comparison
principle, implies the comparison principle for viscosity solutions
immediately. Our arguments for (\ref{u=u}) mainly rely on the
remarkable result Bank and Baum \cite{BB}, which was extended to
nonlinear case in \cite{STZ-duality}.

We also observe that our results make strong use of the representation
of the solution of the semilinear PPDE by means of the corresponding
backward SDEs~\cite{PP1}. This is a serious limitation of our approach
that we hope to overcome in some future work. However, our approach is
suitable for a large class of PPDEs as Hamilton--Jacobi--Bellman
equations, which are related to stochastic control problems, and their
extension to Hamilton--Jacobi--Bellman--Isaacs equations corresponding
to differential games.

The rest of the paper is organized as follows. In Section \ref
{sect-pathwise} we introduce the framework of \cite{Dupire} and \cite
{CF} and adapt it to our problem. We define classical and viscosity
solutions of PPDE in Section \ref{sect-defn}. In Section \ref
{sect-main} we introduce the main results, and in Section \ref
{sect-proof} we prove some basic properties of the solutions, including
existence, stability and the partial comparison principle of viscosity
solutions. Finally in Section \ref{sect-peron} we prove (\ref{u=u})
and the comparison principle for viscosity solutions.

\section{A pathwise stochastic analysis}
\label{sect-pathwise}

In this section we introduce the spaces on which we will define the
solutions of path dependent PDEs. The key notions of derivatives were
proposed by Dupire \cite{Dupire} who introduced the functional It\^o
calculus, and further developed by Cont and Fournie \cite{CF}. We
shall also introduce their localization version for our purpose.

\subsection{Derivatives on c\`{a}dl\`{a}g paths}

Let $\hat\O:= \dbD([0,T], \dbR^d)$, the set of c\`{a}dl\`{a}g
paths, $\hat
\omega$ denote the elements of $\hat\O$, $\hat B$ the canonical\vadjust{\goodbreak}
process, $\hat\dbF$ the filtration generated by $\hat B$ and $\hat\L
:= [0,T]\times\hat\O$.
We define seminorms on $\hat\O$ and a pseudometric on $\hat\L$ as follows:
for any $(t, \hat\omega), ( t', \hat\omega') \in\hat\L$,
%
\begin{eqnarray}
\label{rho}
\|\hat\omega\|_{t}&:=& \sup_{0\le s\le t} |\hat
\omega_s|,\nonumber\\[-8pt]\\[-8pt]
d_\infty\bigl((t, \hat\omega),\bigl(
t', \hat\omega'\bigr) \bigr)&:=& \bigl|t-t'\bigr| +
\sup_{0\le s\le T} \bigl|\hat\omega_{t\wedge s} - \hat\omega'_{t'\wedge
s}\bigr|.\nonumber
\end{eqnarray}
Then $(\hat\O, \|\cdot\|_{T})$ is a Banach space and $(\hat\L, d_\infty)$ is
a complete pseudometric space.

Let $\hat u\dvtx  \hat\L\to\dbR$ be an $\hat\dbF$-progressively
measurable random field. Note that the progressive measurability
implies that $\hat u(t,\hat\omega) = \hat u(t,\hat\omega_{\cdot
\wedge t})$ for all $(t,\hat\omega) \in\hat\L$. Following Dupire
\cite{Dupire}, we define spatial derivatives of $\hat u$, if they
exist, in the standard sense: for the basis $e_i$ of $\dbR^d$,
$i=1,\ldots, d$,
%
\begin{eqnarray}
\label{hatpax} \partial_{\omega_i}\hat u(t,\hat\omega)&:=&
\lim_{h\to0}\frac
{1}{h} \bigl[ \hat u(t,\hat\omega+ h
\one_{[t,T]}e_i) - \hat u(t,\hat\omega) \bigr],\nonumber\\[-8pt]\\[-8pt]
\partial_{\omega_i\omega_j}\hat u&:=& \partial_{\omega
_i}(\hat u_{\omega_j}),\qquad
i,j=1,\ldots,d,\nonumber
\end{eqnarray}
and the right time-derivative of $\hat u$, if it exists, as
%
\begin{equation}
\label{hatpat} \partial_t \hat u(t,\hat\omega):=
\lim_{h\to0, h>0} \frac
{1}{h} \bigl[\hat u (t+h, \hat
\omega_{\cdot\wedge t} )-\hat u (t, \hat\omega) \bigr],\qquad t<T.
\end{equation}
For the final time $T$, we define
%
\begin{equation}
\partial_t \hat u(T,\omega) := \lim_{t<T, t\uparrow T}
\partial_t \hat u(t,\omega).
\end{equation}
We take the convention that $\hat\omega$ are column vectors, but
$\partial_{\omega}\hat u$ denotes row vectors, and $\partial^2_{\omega
\omega} \hat u$ denote $d\times d$-matrices.
%
\begin{defn}
\label{defn-hatspace} Let $\hat u\dvtx  \hat\L\to\dbR$ be $\hat\dbF
$-progressively measurable.

\begin{longlist}[(iii)]
\item[(i)]
We say $\hat u \in C^0( \hat\L)$ if $\hat u$ is continuous in
$(t,\hat\omega)$ under $d_\infty$.

\item[(ii)] We say $\hat u \in C^0_b( \hat\L)\subset C^0(\hat\L)$ if $\hat
u$ is bounded.

\item[(iii)] We say $\hat u\in C^{1,2}_b( \hat\L) \subset C^0( \hat\L)$
if $\partial_t \hat u$, $\partial_{\omega} \hat u$, and $\partial
^2_{\omega
\omega}\hat u$ exist and are in $C^0_b(\hat\L)$.
\end{longlist}
\end{defn}
%
\begin{rem}
\label{rem-bound}
To simplify the presentation, in this paper we will consider only
bounded viscosity solutions. By slightly more involved estimates, we
can extend our results to the cases with polynomial growth. However,
the boundedness of the derivatives $\partial_t \hat u$, $\partial
_{\omega}
\hat u$, and $\partial^2_{\omega\omega}\hat u$ is crucial for the
functional It\^{o}'s formula (\ref{Ito}) below.
\end{rem}
%
\subsection{Derivatives on continuous paths}

We now let $\O:= \{\omega\in C([0,T], \dbR^d)\dvtx\allowbreak \omega_0={\mathbf
0} \}$, the set of continuous paths with initial value ${\mathbf0}$,
$B$ the canonical process, $\dbF$ the filtration generated by
$B$,\vadjust{\goodbreak}
$\dbP_0$ the Wiener measure, and $\L:= [0,T]\times\O$. Here and in
the sequel, for notational simplicity, we use ${\mathbf0}$ to denote
vectors or matrices with appropriate dimensions whose components are
all equal to $0$.

Clearly $\O\subset\hat\O$, $\L\subset\hat\L$, and each $\omega
\in\O$ can also be viewed as an element of~$\hat\O$. Then \mbox{$\|\cdot
\|_t$} and $d_\infty$ in (\ref{rho}) are well defined on $\O$ and $\L
$, $(\O, \|\cdot\|_T)$ is a Banach space, and $(\L, d_\infty)$ is a complete
pseudometric space. Given $u\dvtx  \L\to\dbR$ and $\hat u\dvtx  \hat\L\to\dbR$,
we say $\hat u$ is consistent with $u$ on $\L$ if
%
\begin{equation}
\label{hatu=u} \hat u(t,\omega) = u(t,\omega)\qquad \mbox{for all }(t,\omega)
\in\L.
\end{equation}

\begin{defn}
\label{defn-space} Let $u\dvtx  \L\to\dbR$ be $\dbF$-progressively measurable.

\begin{longlist}[(iii)]
\item[(i)]
We say $u \in C^0( \L)$ if $u$ is continuous in $(t,\omega)$
under $d_\infty$.

\item[(ii)] We say $u \in C^0_b( \L)\subset C^0(\L)$ if $u$ is bounded.

\item[(iii)] We say $u \in C^{1,2}_b(\L)$ if there exists $\hat u \in
C^{1,2}_b( \hat\L)$ such that (\ref{hatu=u}) holds.
\end{longlist}
\end{defn}

By \cite{Dupire} and \cite{CF}, we have the following important results.
%
\begin{theorem}
\label{thm-Ito}
Let $u\in C^{1,2}_b(\L)$ and $\hat u\in C^{1,2}_b( \hat\L)$ such
that (\ref{hatu=u}) holds.

\begin{longlist}[(ii)]
\item[(i)]
The following definition
\[
\partial_t u:= \partial_t \hat u,\qquad \partial_{\omega}
u:= \partial_{\omega} \hat u,\qquad \partial^2_{\omega\omega}u:=
\partial^2_{\omega\omega} \hat u\qquad \mbox{on } \L
\]
is independent of the choice of $\hat u$. Namely, if there is another
$\hat u'\in C^{1,2}_b( \hat\L)$ such that (\ref{hatu=u}) holds, then
the derivatives of $\hat u'$ coincide with those of $\hat u$ on~$\L$.

\item[(ii)] If $\dbP$ is a semimartingale measure, then $u$ is a
semimartingale under $\dbP$ and
%
\begin{equation}
\label{Ito} d u_t= \partial_t u_t \,dt +
\tfrac12\operatorname{tr} \bigl(\partial^2_{\omega\omega}
u_t \,d\langle B\rangle_t \bigr)+ \partial_{\omega}u_t\,dB_t,\qquad
\dbP\mbox{-a.s.}
\end{equation}
\end{longlist}
\end{theorem}
We note that, for any given $\dbP$, the quadratic variation $\langle
B\rangle
$ is well defined. In fact, although not used in this paper, one can
construct $\langle B\rangle$ in a pathwise manner, see, for example, \cite{STZ-2BSDE}.
Here and in the sequel, when we emphasize that $u$ is a process, we use
the notation $u_t(\omega):= u(t,\omega)$ and often omit $\omega$
by simply writing it as $u_t$. Moreover, when a probability is
involved, quite often we use $B$ which by definition satisfies
$B_t(\omega) = \omega_t$.

\subsection{Localization of the spaces}

For our purpose, we need to introduce the localization version of the
above notions. Let
%
\begin{eqnarray}
\label{cT} \cT&:=& \bigl\{\mbox{$\dbF$-stopping time $\t$: for all }t
\in[0,T),
\nonumber\\[-8pt]\\[-8pt]
&&\hspace*{5.2pt}\bigl\{\omega\dvtx  \t(\omega) > t\bigr\} \mbox{ is an open subset of } \bigl(\O,
\|\cdot
\|_T\bigr) \bigr\}.
\nonumber
\end{eqnarray}
The following is a typical example of such $\t$.
%
\begin{eg}
\label{eg-tau}
Let $u\in C^0(\L)$. Then, for any constant $c$,
\[
\t:= \inf\bigl\{t\dvtx  u(t,\omega) \ge c \bigr\}\wedge T\mbox{ is in }\cT.
\]
\end{eg}
\begin{pf}
For any $t<T$, $\{\t> t\} = \{ \sup_{0\le s\le t} u_s < c\}$. Fix
$\omega\in\{\t>t\}$, and set $\e:= \frac12[c-\sup_{0\le s\le t}
u(s,\omega)]>0$. For any $s\in[0,t]$, since $u$ is continuous at
$(s,\omega)$, there exists a constant $h_s>0$ such that $|u(r,\tilde
\omega) - u(s,\omega)|\le\e$ whenever $d_\infty((r,\tilde
\omega), (s,\omega))< h_s$. Note that the open intervals $(s-\frac
12 h_s, s+\frac12 h_s)$, $s\in[0,t]$, cover the compact set $[0,t]$.
Then there exist $0=s_0<s_1<\cdots<s_n=t$ such that $[0,t] \subset
\bigcup_{0\le i\le n} (s_i-\frac12 h_{s_i}, s_i+\frac12 h_{s_i})$. Now set
$h:= \frac12 \min_{0\le i\le n} h_{s_i} > 0$. For any $\tilde
\omega\in\O$ such that $\|\tilde\omega-\omega\|_T <h$, for any
$s\in[0,t]$, there exists $i$ such that $|s-s_i|\le\frac12 h_{s_i}$. Then
\[
d_\infty\bigl((s,\tilde\omega), (s_i, \omega)\bigr)
\le|s-s_i| + \| \tilde\omega-\omega\|_T \le\tfrac12
h_{s_i} + h \le h_{s_i} \qquad\mbox{for all } s\in[0,t].
\]
Thus
\[
u(s,\tilde\omega) \le u(s_i,\omega) + \e\le\sup_{0\le s\le t}
u(s,\omega) + \e< c \qquad\mbox{for all } s\in[0,t].
\]
This implies that $\t(\tilde\omega) > t$, and therefore $\t\in\cT$.
\end{pf}

Denote
%
\begin{equation}
\label{Ltau}\hspace*{28pt} \L(\t):= \bigl\{ (t, \omega)\in\L\dvtx  t<\t(\omega) \bigr\}
\quad\mbox{and}\quad
\bar\L(\t):= \bigl\{ (t, \omega)\in\L\dvtx  t\le\t(\omega) \bigr\}.
\end{equation}
Then clearly $\L(\t)$ is an open subset of $(\L, d_\infty)$.
%
\begin{defn}
\label{defn-local}
Let $\t\in\cT$ and $u\dvtx  \bar\L(\t) \to\dbR$. We say $u \in
C^{1,2}_b(\bar\L(\t))$ if there exists $\tilde u \in C^{1,2}_b(\L)$
such that
%
\begin{equation}
\label{u=tildeu} u = \tilde u \qquad\mbox{on } \bar\L(\t).
\end{equation}
\end{defn}

The following result is the localization version of Theorem \ref{thm-Ito}.
%
\begin{prop}
\label{prop-local}
Let $\t\in\cT$, $u \in C^{1,2}_b(\bar\L(\t))$, $\tilde u\in
C^{1,2}_b(\L)$ such that (\ref{u=tildeu}) holds.

\begin{longlist}[(iii)]
\item[(i)]
One may define
%
\begin{equation}
\label{derivative-local} \partial_tu:= \partial_t\tilde
u,\qquad \partial_{\omega} u:= \partial_{\omega} \tilde u,\qquad
\partial^2_{\omega\omega} u:= \partial^2_{\omega
\omega}
\tilde u\qquad\mbox{on }\L(\t),
\end{equation}
and the definition is independent of the choice of $\tilde u$.

\item[(ii)] Let $\dbP$ be a semimartingale measure. Then $u$ is a $\dbP
$-semimartingale on $[0,\t]$ and (\ref{Ito}) holds on $[0,\t]$.
\end{longlist}
\end{prop}
\begin{pf} First, for the derivatives defined in (\ref
{derivative-local}), (\ref{Ito}) follows directly from Theorem \ref
{thm-Ito}. Next, assume $\tilde u' \in C^{1,2}_b(\L)$ also satisfies
(\ref{u=tildeu}). Denote $\bar u:= \tilde u - \tilde u'$. Then $\bar
u=0$ on $\bar\L(\t)$. Now fix $(t,\omega) \in\L(\t)$. Since $\L
(\t)$ is open, there exists $h:=h(t,\omega)>0$ such that $(s,\tilde
\omega) \in\L(\t)$ whenever $d_\infty((s,\tilde\omega), (t,
\omega)) <h$. Now following the definition of the time derivative we
obtain immediately that $\partial_t \bar u(t,\omega) = 0$. Moreover, let
$\dbP=\dbP_0$, and applying (\ref{Ito}) to $\bar u$, we have
\[
0 = \tfrac12\operatorname{tr} \bigl(\partial^2_{\omega\omega} \bar
u_t \bigr) \,dt+ \partial_{\omega}\bar u_t\,dB_t,\qquad 0
\le t< \t, \dbP_0\mbox{-a.s.}
\]
Thus, since $\partial_{\omega} \bar u$ and $\partial^2_{\omega
\omega}\bar
u$ are bounded,
\[
\partial_{\omega} \bar u = {\mathbf0},\qquad \partial^2_{\omega\omega
}
\bar u = {\mathbf0},\qquad dt\times d\dbP_0\mbox{-a.s. on }\L(\t).
\]
Since $\L(\t)$ is open, and $\partial_{\omega} \bar u$ and
$\partial^2_{\omega\omega} \bar u$ are continuous in $(t,\omega)$ under
$d_\infty$, it is clear that
\[
\partial_{\omega} \bar u = {\mathbf0},\qquad \partial^2_{\omega\omega
}
\bar u = {\mathbf0}\qquad\mbox{on }\L(\t).
\]
This implies that the definition in (\ref{derivative-local}) is
independent of the choice of $\tilde u$.
\end{pf}

\subsection{Space shift}

We first fix $t\in[0,T]$ and introduce the shifted spaces on c\`{a}dl\`{a}g paths:

- Let $\hat\O^t:= \dbD([t,T], \dbR^d)$ be the shifted canonical
space; $\hat B^{t}$ the shifted canonical process on
$\hat\O^t$; $\hat\dbF^{t}$ the shifted filtration generated by
$B^{t}$; and $\hat\L^t:= [t,T]\times\hat\O^t$.

- Define $\|\cdot\|^t_s$ and $d^t_\infty$ in the spirit of (\ref{rho}).

- For $\hat\dbF^t$-progressively measurable $\hat u\dvtx  \hat\L^t \to
\dbR$, define the derivatives in the spirit of (\ref{hatpax}) and
(\ref{hatpat}), and define the spaces $C^0(\hat\L^t)$, $C_b^0(\hat
\L^t)$ and $C^{1,2}_b(\hat\L^t)$ in the spirit of Definition \ref
{defn-hatspace}.\vspace*{9pt}



Similarly, we may define the shifted spaces on continuous
paths:\vspace*{9pt}

- Let $\O^t:= \{\omega\in C([t,T], \dbR^d)\dvtx  \omega_t ={\mathbf
0} \}$ be the shifted canonical space, $B^{t}$ the shifted
canonical process on
$\O^t$, $\dbF^{t}$ the shifted filtration generated by $B^{t}$, $\dbP
^t_0$ the Wiener measure on $\O^t$ and $\L^t:= [t,T]\times\O^t$.

- Define $C^0(\L^t)$, $C_b^0(\L^t)$ and $C^{1,2}_b(\L^t)$ in an
obvious way.

- Let $\cT^t$ denote the space of $\dbF^t$-stopping times $ \t$ such
that, for any $s\in[t,T)$, the set $\{\omega\in\O^t\dvtx  \t(\omega)
> s\}$ is an open subset of $\O^t$ under $\|\cdot\|^t_T$.

- For each $\t\in\cT^t$, define $\L^t(\t)$, $\bar\L^t(\t)$, and
$C^{1,2}_b(\bar\L^t(\t))$ in an obvious way.\vspace*{9pt}

We next introduce the shift and concatenation operators. Let $0\le
s\le\break
t\le T$.\vspace*{9pt}


- For $\hat\omega\in\hat\O^s$, $\hat\omega' \in\hat\O^t$
and $\omega\in\O^s$, $\omega'\in\O^t$, define the concatenation
paths $\hat\omega\otimes_{t} \hat\omega'\in\hat\O^s$ and
$\omega\otimes_{t} \omega'\in\O^s$ by
\begin{eqnarray*}
\bigl(\hat\omega\otimes_t \hat
\omega'\bigr) (r)&:=&\hat\omega_r\one_{[s,t)}(r)
+ \bigl(\hat\omega_{t-} 
+ \hat\omega'_r\bigr)\one_{[t, T]}(r);
\\[-10pt]
&&\hspace*{148pt}\qquad\mbox{for all } r\in[s,T].
\\[-10pt]
\bigl(\omega\otimes_t \omega'
\bigr) (r)&:=& \omega_r\one_{[s,t)}(r) + \bigl(
\omega_{t} + \omega'_r\bigr)
\one_{[t, T]}(r);
\end{eqnarray*}

- Let $\hat\omega\in\hat\O^s$. For $\hat\cF^{s}_{T}$-measurable
random variable $\hat\xi$ and $\hat\dbF^{s}$-progressively measurable
process $\hat X$ on $\hat\O^s$, define the shifted $\hat\cF
^{t}_{T}$-measurable random variable $\hat\xi^{t, \hat\omega}$ and
$\hat\dbF^{t}$-progressively measurable
process $\hat X^{t,\hat\omega}$ on $\hat\O^t$ by
\[
\hat\xi^{t, \hat\omega}\bigl(\hat\omega'\bigr):=\hat\xi\bigl(\hat
\omega\otimes_t \hat\omega'
\bigr),\qquad \hat X^{t, \hat\omega}\bigl(\hat\omega'\bigr):=\hat X\bigl(
\hat\omega\otimes_t \hat\omega'
\bigr)\qquad\mbox{for all } \hat\omega'\in\hat\O^t.
\]

- Let $\omega\in\O^s$. For $\cF^{s}_{T}$-measurable
random variable $\xi$ and $\dbF^{s}$-progressively measurable
process $X$ on $\O^s$, define the shifted $\cF^{t}_{T}$-measurable
random variable $\xi^{t, \omega}$ and $\dbF^{t}$-progressively measurable
process $X^{t,\omega}$ on $\O^t$ by
\[
\xi^{t, \omega}\bigl(\omega'\bigr):=\xi\bigl(\omega\otimes
_t \omega'\bigr),\qquad X^{t, \omega}
\bigl(\omega'\bigr):= X\bigl(\omega\otimes
_t \omega'\bigr)\qquad\mbox{for all }
\omega'\in\O^t.
\]

It is clear that all the results in previous subsections can be
extended to the shifted spaces, after obvious modifications. Moreover,
for any $\t\in\cT$, $(t, \omega) \in\L(\t)$ and $u\in
C^{1,2}_b(\bar\L(\t))$, it is clear that $\t^{t,\omega} \in\cT^t$ and
$u^{t,\omega} \in C^{1,2}_b(\bar\L^t(\t^{t,\omega}))$.

For some technical proofs later, we shall also use the following space.
Denote
%
\begin{equation}
\label{Tt+} \cT^t_+:= \bigl\{\t\in\cT^t\dvtx  \t>t\bigr\}
\qquad\mbox{for }t<T \mbox{ and } \cT^T_+:= \{T\}.
\end{equation}



%
\begin{defn}
\label{defn-barCP}
Let $t\in[0,T]$, $u\dvtx  \L^t \to\dbR$ and $\dbP$ be a semimartingale
measure on $\O^t$. We say $u\in\bar C^{1,2}_\dbP(\L^t)$ if there
exist an increasing sequence of $\dbF^t$-stopping times $t=\t_0 \le
\t_1 \le\cdots\le T$ such that:

\begin{longlist}[(iii)]
\item[(i)]
For each $i\ge0$ and $\omega\in\O^t$,
\[
\t^{\t_i(\omega
),\omega}_{i+1} \in\cT^{\t_i(\omega)}_+\quad \mbox{and}\quad u^{\t_i(\omega
),\omega} \in C^{1,2}_b\bigl(\bar\L^{\t_i(\omega)}\bigl( \t^{\t
_i(\omega),\omega}_{i+1}\bigr)\bigr);
\]

\item[(ii)] For each $i\ge0$ and $\omega\in\O$, $u_\cdot(\omega)$ is
continuous on $[0,\t_i(\omega)]$;

\item[(iii)] For $\dbP$-a.s. $\omega\in\O^t$, the set $\{i\dvtx  \t_i(\omega
) < T\}$ is finite.
\end{longlist}
\end{defn}
We shall emphasize that, for $u\in\bar C^{1,2}_\dbP(\L^t)$, the
derivatives of $u$ are bounded on each interval $[\t_i(\omega), \t^{\t
_i(\omega),\omega}_{i+1}]$; however, in general they may be
unbounded on the whole interval $[t,T]$. Also, the previous definition
and, more specifically the dependence on $\dbP$ introduced in item
(iii), is motivated by the results established in Section \ref
{sect-peron} below.

The following result is a direct consequence of Proposition \ref{prop-local}.
%
\begin{prop}
\label{prop-local-Ito}
Let $\dbP$ be a semimartingale measure on $\O^t$ and $u\in\bar
C^{1,2}_\dbP(\L^t)$. Then $u$ is a local $\dbP$-semimartingale on
$[t,T]$ and
\[
d u_s= \partial_t u_s \,ds+
\tfrac12\operatorname{tr} \bigl(\partial^2_{\omega\omega}
u_s\,d\bigl\langle B^t\bigr\rangle_s \bigr)+
\partial_{\omega} u_s\,dB^t_s,\qquad t\le s\le
T, \dbP\mbox{-a.s.}
\]
\end{prop}

\section{PPDEs and definitions}
\label{sect-defn}

In this paper we study the following semi-linear parabolic
Path-dependent PDE (PPDE, for short):
%
\begin{eqnarray}
\label{PPDE} &&(\cL u) (t,\omega) =0,\qquad 0\le t< T,\omega\in\O;
\nonumber\\
&&\qquad\mbox{where } (\cL u) (t,\omega):= -\partial_t u(t,\omega) - \tfrac
12 \operatorname{tr} \bigl(\partial^2_{\omega\omega} u(t,\omega)
\bigr) \\
&&\qquad\hspace*{3pt}\hphantom{\mbox{where } (\cL u) (t,\omega):=}{}- f\bigl(t,\omega, u(t,\omega), \partial_{\omega} u(t,\omega)
\bigr).
\nonumber
\end{eqnarray}
We remark that there is a potential to extend our results to a much
more general setting. 
However, in order to focus on the main ideas, in this paper we content
ourselves with the simple PPDE (\ref{PPDE}) under somewhat strong
technical conditions, and leave more general cases, for example, fully
nonlinear PPDEs, for future studies.
%
\begin{rem}
\label{rem-markovian}
In the Markovian case, namely $f = f(t, \omega_t, y, z)$ and
$u(t,\allowbreak\omega) = v(t,\omega_t)$, the PPDE (\ref{PPDE}) reduces to the
following PDE:
%
\begin{eqnarray}
\label{PDE} &&(\cL v) (t,x) =0,\qquad 0\le t< T,x\in\dbR^d,
\nonumber\\
&&\qquad\mbox{where } (\cL v) (t,x):= -\partial_t v(t,x) - \tfrac12
\operatorname{tr} \bigl[D^2_{xx} v(t,x) \bigr]\\
&&\qquad\hspace*{3pt}\hphantom{\mbox{where } (\cL v) (t,x):=}{} - f
\bigl(t,x, v(t,x), D_x v(t,x)\bigr).
\nonumber
\end{eqnarray}
Here $D_x$ and $D^2_{xx}$ denote the standard first and second order
derivatives with respect to $x$. However, slightly different from the
PDE literature but consistent with~(\ref{hatpat}), $\partial_t$
denotes the
right time-derivative.
\end{rem}

As usual, we start with classical solutions.
%
\begin{defn}
\label{defn-classical} Let $u\in C^{1,2}_b(\L)$. We say $u$ is a
classical solution (resp., sub-solution, super-solution) of PPDE (\ref
{PPDE}) if
%
\begin{equation}
\label{classical-subsuper} (\cL u) (t,\omega) = (\mbox{resp.,}\le, \ge)
0
\qquad\mbox{for all } (t,\omega) \in[0,T)\times\O.
\end{equation}
\end{defn}
It is clear that, in the Markovian setting as in Remark \ref{rem-markovian},
\begin{eqnarray*}
\begin{tabular}{p{330pt}}
$u$ is a classical solution (resp., sub-solution,
super-solution) of
PPDE (\ref{PPDE})
if and only if $v$ is a classical solution (resp., sub-solution,
super-solution) of PDE~(\ref{PDE}).
\end{tabular}
\end{eqnarray*}

Existence and uniqueness of classical solutions are related to the
analogue results for the corresponding backward SDE. In order to avoid
diverting the attention from our main purpose in this paper, we report
these properties later in Section \ref{sect-classical}, and we move
to our notion of viscosity solutions.

For any $L\ge0$ and $t<T$, let $\cU^L_t$ denote the space of $\dbF
^t$-progressively measurable $\dbR^d$-valued processes $\b$ such that
each component of $\b$ is bounded by $L$. By viewing $\b$ as row
vectors, we define
%
\begin{eqnarray}
\label{Mbeta} M^{t,\b}_s:= \exp\biggl(\int
_t^s \b_r \,dB_r^t
-\frac12 \int_t^s |\b_r|^2\,dr
\biggr),\nonumber\\[-8pt]\\[-8pt]
&&\eqntext{\dbP_0^t\mbox{-a.s.}, d\dbP^{t,\b}:=
M^{t,\b}_T \,d\dbP_0^t,}
\end{eqnarray}
and we introduce for all $t\in[0,T]$ two nonlinear expectations: for
any $\xi\in L^2(\cF^t_T, \dbP^t_0)$,
%
\begin{eqnarray}\label{cE}
\underline\cE^L_t[\xi] &:=& \inf\bigl\{
\dbE^{\dbP^{t,\b}}[\xi]\dvtx  \b\in\cU^L_t \bigr\};
\nonumber\\[-8pt]\\[-8pt]
\overline \cE{}^L_t[\xi] &:=& \sup\bigl\{\dbE^{\dbP^{t,\b}}[\xi]\dvtx  \b
\in\cU^L_t \bigr\}.\nonumber
\end{eqnarray}
Moreover, for any $u\in C^0_b(\L)$, define
%
\begin{eqnarray}
\label{cA} &&\underline\cA^{L}u (t,\omega) := \Bigl\{\f\in
C^{1,2}_b\bigl(\L^t\bigr)\mbox{: there exists $
\t\in\cT^t_+$ such that}
\nonumber
\\
&&\hspace*{68pt}0=\f(t,{\mathbf0})-u(t,\omega) =\min_{\tilde\t\in\cT^t} \underline\cE^{L}_t
\bigl[\bigl(\f- u^{t,\omega}\bigr)_{\tilde\t\wedge\t} \bigr] \Bigr\};
\nonumber\\[-8pt]\\[-8pt]
&&\overline\cA{}^{L}u (t,\omega) := \Bigl\{\f\in C^{1,2}_b
\bigl(\L^t\bigr)\mbox{: there exists $\t\in\cT^t_+$ such
that}
\nonumber
\\
&&\hspace*{68pt}0=\f(t,{\mathbf0})-u(t,\omega) =\max_{\tilde\t\in\cT^t} \overline\cE{}^{L}_t
\bigl[\bigl(\f- u^{t,\omega}\bigr)_{\tilde\t\wedge\t} \bigr]
\Bigr\}.\nonumber
\end{eqnarray}

\begin{defn}
\label{defn-viscosity}
Let $u\in C_b^0(\L)$.

\begin{longlist}[(iii)]
\item[(i)]
For any $L\ge0$, we say $u$ is a viscosity $L$-subsolution
(resp.,\vspace*{1pt}
$L$-supersolu\-tion) of PPDE (\ref{PPDE}) if, for any $(t,\omega)\in
[0, T)\times\O$ and any $\f\in\underline\cA^{L}u (t,\omega)$
[resp., $\f\in\overline\cA{}^{L}u (t,\omega)$], it holds that
\[
\bigl(\cL^{t,\omega}\f\bigr) (t,{\mathbf0}) \le (\mbox{resp.,} \ge)0,
\]
where, for each $(s,\tilde\omega) \in[t,T]\times\O^t$,
\[
\bigl(\cL^{t,\omega} \f\bigr) (s,\tilde\omega):=- \partial_t
\f(s,\tilde\omega) - \tfrac12 \operatorname{tr} \bigl[ \partial
^2_{\omega
\omega}
\f(s,\tilde\omega) \bigr] - f^{t,\omega} \bigl(s,\tilde\omega, \f
(s,\tilde
\omega), \partial_{\omega} \f(s,\tilde\omega) \bigr).
\]

\item[(ii)] We say $u$ is a viscosity subsolution (resp., supersolution) of
PPDE (\ref{PPDE}) if $u$ is viscosity $L$-subsolution (resp.,
$L$-supersolution) of PPDE (\ref{PPDE}) for some $L\ge0$.

\item[(iii)] We say $u$ is a viscosity solution of PPDE (\ref{PPDE}) if it is
both a viscosity subsolution and a viscosity supersolution.
\end{longlist}
\end{defn}

In the rest of this section we provide several remarks concerning our
definition of viscosity solutions. In most places we will comment on
the viscosity subsolution only, but obviously similar properties hold
for the viscosity supersolution as well.
%
\begin{rem}
\label{rem-cA}
As standard in the literature on viscosity solutions of PDEs:

\begin{longlist}[(iii)]
\item[(i)] The viscosity property is a local property in the following sense.
For any $(t,\omega) \in[0,T)\times\O$ and any $\e>0$, define
\[
\t_\e:= \inf\bigl\{ s>t\dvtx  \bigl|B^t_s\bigr| \ge\e
\bigr\} \wedge(t+\e).
\]
To check the viscosity property of $u$ at $(t,\omega)$, it suffices
to know the value of $u^{t,\omega}$ on $[t, \t_\e]$ for an
arbitrarily small $\e>0$.

\item[(ii)] Typically $\underline\cA^{L}u (t,\omega)$ and $\overline\cA{}^{L}u
(t,\omega)$ are disjoint, so $u$ is a viscosity solution does
not mean $(\cL^{t,\omega}\f)(t,{\mathbf0}) = 0$ for $\f$ in some
appropriate set. One has to check viscosity subsolution property and
viscosity supersolution property separately.

\item[(iii)] In general $\underline\cA^{L}u (t,\omega)$ could be empty. In
this case automatically $u$ satisfies the viscosity subsolution
property at $(t,\omega)$.
\end{longlist}
\end{rem}
%
\begin{rem}
\label{rem-L}
(i) For $0\le L_1 < L_2$, obviously $\cU^{L_1}_t \subset\cU
^{L_2}_t$, $\underline\cE^{L_2}_t\le\underline\cE^{L_1}_t$ and
$\underline\cA^{L_2}u(t,\omega) \subset\underline\cA^{L_1}u(t,\omega)$.
Then one can easily check that a viscosity
$L_1$-subsolution must be a viscosity $L_2$-subsolution. Consequently,
\begin{eqnarray*}
\begin{tabular}{p{330pt}}
$u$ is a viscosity subsolution
if and only if there exists a $L\ge0$ such that, for all
$\tilde
L\ge L$, $u$ is a viscosity $\tilde L$-subsolution.
\end{tabular}
\end{eqnarray*}

(ii) However, we require the same $L$ for all $(t,\omega)$. We should
point out that our definition of viscosity subsolution is not
equivalent to the following alternative definition, under which we are
not able to prove the comparison principle:
\[
\mbox{for any $(t,\omega)$ and any $\f\in\bigcap_{L\ge0}
\underline\cA^Lu(t,\omega)$, it holds that $(\cL^{t,\omega}
\f) (t,{\mathbf0}) \le0$.}
\]
\end{rem}
%
\begin{rem}
\label{rem-strict}
We may replace $\underline\cA^{L}$ with the following
$(\underline\cA')^{L}$ which requires strict inequality,
%
\begin{eqnarray}
\label{cA1}\qquad
&&{\underline\cA'}^{L}u (t,
\omega) := \bigl\{\f\in C^{1,2}_b\bigl(\L^t
\bigr)\mbox{: there exists $\t\in\cT^t_+$ such that}
\nonumber\\[-8pt]\\[-8pt]
&&\hspace*{69pt}0=\f(t,{\mathbf0})-u(t,\omega) <\underline\cE^{L}_t \bigl[
\bigl(\f- u^{t,\omega}\bigr)_{\tilde\t\wedge\t} \bigr] \mbox{ for
all }\tilde\t\in
\cT^t_+ \bigr\}.\nonumber
\end{eqnarray}
Then $u$ is a viscosity $L$-subsolution of PPDE (\ref{PPDE}) if and
only if
\[
\mbox{$\bigl(\cL^{t,\omega}\f\bigr) (t,{\mathbf0}) \le0$}\qquad\mbox{for all $(t,\omega
)
\in[0,T)\times\O$ and $\f\in{\underline\cA'}^{L}u (t,
\omega)$.}
\]
A similar statement holds for the viscosity supersolution.

Indeed, since ${\underline\cA'}^{L} u(t,\omega) \subset\underline
\cA^{L}u (t,\omega)$, then only the if part is clear. To prove the
if part, let $(t,\omega) \in[0,T)\times\O$ and $\f\in\underline
\cA^{L}u (t,\omega)$. For any $\e>0$, denote $\f^\e(s, \tilde
\omega):= \f(s,\tilde\omega) + \e(s-t)$. Then clearly $\f^\e
\in{\underline\cA'}^{L}u (t,\omega)$, and thus
\begin{eqnarray*}
\bigl(\cL^{t,\omega} \f^\e\bigr) (t,{\mathbf0}) &=& -
\partial_t \f(t,{\mathbf0}) -\e- \tfrac12 \operatorname{tr} \bigl(
\partial^2_{\omega\omega} \f(t,{\mathbf0}) \bigr) \\
&&{}- f^{t,\omega}
\bigl(t,\omega, \f(t,{\mathbf0}), \partial_{\omega} \f(t,{\mathbf0}) \bigr) \\
&\le&0.
\end{eqnarray*}
Send $\e\to0$, we obtain $(\cL^{t,\omega} \f)(t,{\mathbf0}) \le0$,
and thus $u$ is a viscosity $L$-subsolution.
\end{rem}
%
\begin{rem}
\label{rem-markovian-viscosity}
Consider the Markovian setting in Remark \ref{rem-markovian}. One
can easily check that $u$ is a viscosity subsolution of PPDE (\ref
{PPDE}) in the sense of Definition \ref{defn-viscosity} implies that
$v$ is a viscosity subsolution of PDE (\ref{PDE}) in the standard sense.
\end{rem}
%
\begin{rem}
\label{rem-cA2}
We have some flexibility to choose $\underline\cA^Lu(t,\omega
)$ and $\overline\cA{}^Lu(t,\omega)$ in Definition \ref
{defn-viscosity}. In principle, the smaller these sets are, the more
easily we can prove viscosity properties and thus the existence of
viscosity solutions, but the comparison principle and the uniqueness of
viscosity solutions become more difficult.

\begin{longlist}[(iii)]
\item[(i)]
The following ${\underline\cA''}^{L}u(t,\omega)$ is larger than
$\underline\cA^{L}u(t,\omega)$, but all the results in this paper
still hold true if we use ${\underline\cA''}^{L}u(t,\omega)$ [and
the corresponding ${\overline\cA''}{}^{L}u(t,\omega)$],
%
\begin{eqnarray}
\label{cA2}
&&{\underline\cA''}^{L}u (t,
\omega) := \bigl\{\f\in C^{1,2}_b\bigl(\L^t
\bigr)\mbox{: for any $\t\in\cT^t_+$,}
\nonumber\\
&&\hspace*{72pt}0=\f(t,{\mathbf0})- u(t,\omega) \le\underline\cE^{L}_t \bigl[
\bigl(\f- u^{t,\omega}\bigr)_{\tilde\t\wedge
\t} \bigr] \\
&&\hspace*{187pt}\mbox{for some }\tilde\t\in
\cT^t_+ \bigr\}.\nonumber
\end{eqnarray}

\item[(ii)] However, if we use the following smaller alternatives of
$\underline\cA^Lu (t,\omega)$, which do not involve the nonlinear
expectation, 
we are not able to prove the comparison principle and the uniqueness of
viscosity solutions,
\begin{eqnarray*}
&&\underline\cA^{\circ}u (t,\omega) := \bigl\{\f\in C^{1,2}_b
\bigl(\L^t\bigr)\mbox{: there exists $\t\in\cT^t_+$ such
that}
\\
&&\hspace*{65pt} 0=\f(t,{\mathbf0})-u(t,\omega)\le\bigl(\f- u^{t,\omega}\bigr)_{\tilde\t
\wedge\t}
\mbox{ for any $\tilde\t\in\cT^t_+$} \bigr\};
\end{eqnarray*}
or
\begin{eqnarray*}
&&\underline\cA^{\circ\circ}u (t,\omega) := \bigl\{\f\in C^{1,2}_b
\bigl(\L^t\bigr)\mbox{: for all }(s,\tilde\omega)\in(t,T]\times
\O^t,
\\
&&\hspace*{77pt} 0=\f(t,{\mathbf0})-u(t,\omega) \le\bigl(\f- u^{t,\omega}\bigr) (s,\tilde
\omega) \bigr\}.
\end{eqnarray*}
See also Remark \ref{rem-L}(ii).
\end{longlist}
\end{rem}
%
\begin{rem}
\label{rem-change}
(i) Let $u$ be a viscosity subsolution of PPDE (\ref{PPDE}). Then
for any $\lambda\in\dbR$, $\tilde u_t:= e^{\lambda t}u_t$ is a
viscosity subsolution of the following PPDE:
%
\begin{equation}
\label{tildePPDE} \tilde\cL\tilde u := - \partial_t \tilde u -
\tfrac12 \operatorname{tr}\bigl(\partial^2_{\omega\omega} \tilde u
\bigr) -\tilde f(t, \omega, \tilde u, \partial_{\omega}\tilde u) \le0,
\end{equation}
where
\[
\tilde f(t, \omega, y,z) := -\lambda y+e^{\lambda t} f\bigl(t,\omega,
e^{-\lambda t}y, e^{-\lambda
t} z\bigr).
\]
Indeed, assume $u$ is a viscosity $L$-subsolution of PPDE (\ref{PPDE}).
Let $(t,\omega) \in[0,T)\times\O$ and $\tilde\f\in\underline
\cA^{ L} \tilde u(t,\omega)$. For any $\e>0$, denote
\[
\f^\e_s:= e^{-\lambda s}\tilde\f_s +
\e(s-t).
\]
Then, noting that $\tilde\f_t =e^{\lambda t} u(t,\omega)$,
\begin{eqnarray*}
&&
\f^\e_s - u^{t,\omega}_s -
e^{-\lambda t} \bigl(\tilde\f_s - \tilde u^{t,\omega}_s
\bigr) \\
&&\qquad= \bigl(e^{-\lambda s} - e^{-\lambda t} \bigr)\tilde\f_s
+ \bigl(e^{\lambda(s-t)}-1 \bigr) u_s + \e(s-t)
\\
&&\qquad= \bigl(e^{-\lambda s} - e^{-\lambda t} \bigr) (\tilde\f_s -
\tilde\f_t ) + \bigl(e^{\lambda(s-t)}-1 \bigr) (u_s-u_t)
\\
&&\qquad\quad{} + \bigl(e^{-\lambda(s-t)} + e^{\lambda(s-t)} -2 \bigr) u_t +
\e(s-t)
\\
&&\qquad\ge \e(s-t) - C(s-t) \bigl(|\tilde\f_s - \tilde\f_t| +
|u_s-u_t| + (s-t) \bigr).
\end{eqnarray*}
Let $\tilde\t\in\cT^t_+$ be a stopping time corresponding to
$\tilde\f\in\underline\cA^{ L} \tilde u(t,\omega)$, and set
\[
\t_\e:= \tilde\t\wedge\inf\biggl\{s>t\dvtx  |\tilde\f_s-
\tilde\f_t| + |u_s-u_t| + (s-t) \ge
\frac{\e}{C} \biggr\}\wedge T.
\]
Then $\t_\e\in\cT^t_+$, by Example \ref{eg-tau}, and for any $\t
\in\cT^t$ such that $\t\le\t_\e$, it follows from the previous
inequality that
\[
\f^\e_\t- u^{t,\omega}_\t\ge
e^{-\lambda t}\bigl[\tilde\f_\t- \tilde u^{t,\omega}_\t
\bigr].
\]
By the increase and the homogeneity of the operator $\underline{\cE
}^L_t$, together with the fact that $\tilde\f\in\underline\cA^{ L}
\tilde u(t,\omega)$, this implies that
\[
\underline\cE^L_t \bigl[\f^\e_\t-
u^{t,\omega}_\t\bigr] \ge\underline\cE^L_t
\bigl[e^{-\lambda t}\bigl(\tilde\f_\t- \tilde u^{t,\omega}_\t
\bigr) \bigr] = e^{-\lambda t} \underline\cE^L_t \bigl[
\tilde\f_\t- \tilde u^{t,\omega}_\t\bigr] \ge0 =
\f_t^\e-u_t.
\]
This implies that $\f^\e\in\underline\cA^L u(t,\omega)$, then
$\cL^{t,\omega}\f^\e(t,0) \le0$. Send $\e\to0$, and similar
to Remark \ref{rem-strict} we get $\cL^{t,\omega}\f^0(t,0) \le0$,
where $\f^0_s:= e^{-\lambda s} \tilde\f_s$.
Now by straightforward calculation we obtain
\[
- \partial_t \tilde\f(t,0) - \tfrac12 \operatorname{tr} \bigl[
\partial^2_{\omega\omega} \tilde\f(t,0) \bigr] - \tilde f \bigl(t,
\omega, \tilde\f(t,0), \partial_{\omega} \tilde\f(t,0) \bigr) \le0.
\]
That is, $\tilde u$ is a viscosity subsolution of PPDE (\ref{tildePPDE}).

(ii) If we consider more general variable change: $\bar u(t,\omega)
:= \psi(t, u(t,\omega))$, where $\psi\in C^{1,2}([0,T]\times\dbR
)$ such that $\partial_y \psi>0$. Denote by $\bar\psi:=\psi^{-1}$ the
inverse function of $\psi$ with respect to the space variable $y$.
Then one can easily check that $u$ is a classical subsolution of PPDE
(\ref{PPDE}) if and only if $\bar u$ is a classical subsolution of the
following PPDE:
%
\begin{eqnarray}
\label{barPPDE}
&&\bar\cL\bar u := - \partial_t \bar u - \frac12
\operatorname{tr}\bigl(\partial^2_{\omega\omega} \bar u\bigr) -\bar
f(t, \omega, \bar u, \partial_{\omega}\bar u) \le0
\nonumber
\\
&&\qquad\mbox{where }\bar f(t, \omega, y,z) := \frac{1}{\partial_y\bar\psi(t, y)} \biggl[
\partial_t \bar
\psi(t, y) + \frac12 \partial^2_{yy}\bar\psi(t,y)
|z|^2\\
&&\qquad\hphantom{\mbox{where }\bar f(t, \omega, y,z) := \frac{1}{\partial_y\bar\psi(t, y)} \biggl[}
{} + f\bigl(t,\omega, \bar\psi(t, y), \partial_y\bar
\psi(t, y)z\bigr) \biggr].\nonumber
\end{eqnarray}
However, if $u$ is only a viscosity subsolution of PPDE (\ref{PPDE}),
we are not able to prove that $\bar u$ is a viscosity subsolution of
(\ref{barPPDE}). The main difficulty is that the nonlinear expectation
$\underline\cE^L_t$ and the nonlinear function $\psi$ do not
commute. Consequently, given $\bar\f\in\underline\cA^L\bar
u(t,\omega)$, we are not able to construct as in (i) the
corresponding $\f\in\underline\cA^L u(t,\omega)$.
\end{rem}

We conclude this section by connecting the nonlinear expectation
operators to backward SDEs, and providing some tools from optimal
stopping theory which will be used later.
%
\begin{rem}[(Connecting $\underline{\cE}^L$ and $\overline{\cE}{}^L$
to backward SDEs)]
\label{rem-cE}
For readers who are familiar with BSDE literature, by the comparison
principle of BSDEs (see, e.g., El Karoui, Peng and Quenez \cite{EPQ}),
one can easily show that $\underline\cE^{L}_t[\xi]=\underline\cY_t$ and
$\overline\cE{}^{L}_t[\xi]=\overline\cY_t$, where
$(\underline\cY, \underline\cZ)$ and $(\overline\cY, \overline
\cZ)$ are the solution to the following BSDEs, respectively:
%
\begin{eqnarray}
\label{BSDEt}\quad \underline\cY_s &=& \xi- \int_s^T
L|\underline\cZ_r|\,dr - \int_s^T
\underline\cZ_r\,dB^t_r,
\nonumber\\[-8pt]\\[-8pt]
\overline\cY_s &=& \xi+ \int_s^T L |
\overline\cZ_r| \,dr - \int_s^\t
\overline\cZ_r\,dB^t_r,\qquad t\le s\le T,
\dbP^t_0\mbox{-a.s.}\nonumber
\end{eqnarray}
Moreover, this is a special case of the so called $g$-expectation; see
Peng \cite{Pengg}.
\end{rem}
%
\begin{rem}[(Optimal stopping under nonlinear expectation and reflected
backward SDEs)]
\label{rem-cEreflected}
The definition of the set $\underline{\cA}^L$ is closely related to the following optimal
stopping problem under nonlinear expectation
\[
Y_t := \inf_{\tilde\tau\in\cT^t} \underline{\cE}^L_t
[X_{\tilde\tau\wedge\tau} ]
\]
for some stopping time $\tau\in\cT^t_+$ and some adapted bounded
pathwise continuous process $X$. For the ease of presentation here, we provide only heuristic
arguments, and we refer to Section~7 of \cite{NewRef1} for a rigorous argument and
to \cite{NewRef2} for the optimal stopping problem under more general nonlinear
expectations. 

For later use, we provide some key results which can be proved by
following the standard corresponding arguments in the standard optimal
stopping theory, and we observe that the process $Y$ is pathwise
continuous; see (iv) below.

Following the classical arguments in optimal stopping theory, we have:

\begin{longlist}[(iii)]
\item[(i)]
$\underline{\cE}^L_t [Y_{\tilde\tau\wedge\tau} ]\ge
Y_t$ for all $\tilde\tau\in\cT^t$, that is, $Y$ is an $\underline
{\cE}^L$-submartingale.\vadjust{\goodbreak}

\item[(ii)] If $\tau^*\in\cT^t$ is an optimal stopping rule, then
\[
Y_t = \underline{\cE}^L_t [X_{\tau^*\wedge\tau}
] = \inf_{\tilde\tau\in\cT^{\tau^*}} \underline{\cE}^L_t
[X_{\tilde\tau\wedge\tau} ] = \inf_{\tilde\tau\in\cT^{\tau^*}}
\underline{\cE}^L_t
[Y_{\tilde\tau\wedge\tau} ] = \underline{\cE}^L_t
[Y_{\tau^*\wedge\tau} ],
\]
where the last inequality is a consequence of (i), and the third
inequality follows from the fact that $X\le Y$ on one hand, and $\inf
\underline{\cE}^L_t[\cdot]\ge\underline{\cE}^L_t[\inf\cdot]$ on the other
hand. This implies that $Y_{\tau^*}=X_{\tau^*}$ and, by (i), that
$Y_{\cdot\wedge\tau^*}$ is an $\underline{\cE}^L$-martingale.

\item[(iii)] We then define $\tau^1_t:=\inf\{s>t\dvtx Y_t=X_t\}$. Since
$Y_T=X_T$, we have \mbox{$\tau^1_t\le T$}, a.s. Moreover, following the
classical arguments in optimal stopping theory, we see that $ \{
Y_{s\wedge\tau_t^1} \}_{s\ge t}$ is an $\underline{\cE
}^L$-martingale. With this in hand, we conclude that $\tau^1_t$ is an
optimal stopping time, that is, $Y_t=\underline{\cE}^L_t [X_{\tau
_t^1} ]$.

\item[(iv)] For those readers who are familiar with backward stochastic
differential equations, we mention that $Y=\underline\cY^{\circ}$,
where $(\underline\cY^{\circ}, \underline\cZ^{\circ},\underline
{K}^{\circ})$ is the solution to the following reflected BSDEs:
%
\begin{eqnarray}
\label{BSDEtreflected}\qquad && \underline\cY^{\circ}_s =
X_\t- \int_s^\t L|\underline
\cZ^{\circ}_r|\,dr - \int_s^\t
\underline\cZ^{\circ}_r\,dB^t_r - \int
_s^\t d\underline{K}^{\circ}_s,
\\
&& \underline\cY^{\circ}_s\le X_s\quad \mbox{and}\quad
\bigl(\underline\cY^{\circ}_s-X_s\bigr)\,d
\underline{K}^{\circ}_s=0,\qquad s\in[t,T],\dbP^t_0
\mbox{-a.s.};
\end{eqnarray}
see, for example, \cite{EKPPQ}. In particular, it is a well-known
result that the process $Y$ is pathwise continuous.

\item[(v)] Similar results hold for $\sup_{\tilde\tau\in\cT^t}\overline
{\cE}_t[X_{\tilde\tau\wedge\tau}]$.
\end{longlist}
\end{rem}

\section{The main results}
\label{sect-main}

We start with a stability result.
%
\begin{theorem}
\label{thm-stability}
Let $(f^\e,\e>0)$ be a family of coefficients converging uniformly
toward a coefficient $f\in C^0(\Lambda)$ as $\e\to0$. For some
$L>0$, let $u^\e$ be a viscosity $L$-subsolution (resp.,
$L$-supersolution) of PPDE (\ref{PPDE}) with coefficients~$f^\e$, for
all $\e>0$. Assume further that $u^\e$ converges to some $u$,
uniformly in $\Lambda$. Then $u$ is a viscosity $L$-subsolution (resp.,
supersolution) of PPDE (\ref{PPDE}) with coefficient~$f$.
\end{theorem}

The proof of this result is reported in Section \ref{sect-stability}.
For our next results, we shall always use the following standing
assumptions, where $g$ is a terminal condition associated to the PPDE
(\ref{PPDE}).
%
\begin{assum}
\label{assum-fg}
%
(i) $f$ is bounded, $\dbF$-progressively measurable, continuous in
$t$, uniformly continuous in $\omega$, and uniformly Lipschitz
continuous in $(y,z)$ with a Lipschitz constant $L_0>0$.

(ii) $g$ is bounded and uniformly continuous in $\omega$.
\end{assum}

To establish an existence result of viscosity solutions under the above
assumption, we note that the PPDE (\ref{PPDE}) with terminal
condition\vadjust{\goodbreak}
$u(T,\omega) = g(\omega)$ is closely related to (and actually
motivated from) the following BSDE:
%
\begin{eqnarray}
\label{BSDE} Y^0_t &=& g(B_\cdot) + \int
_t^T f\bigl(s, B_\cdot,
Y^0_s, Z^0_s\bigr) \,ds \nonumber\\[-8pt]\\[-8pt]
&&{}- \int
_t^T Z^0_s
\,dB_s,\qquad 0\le t\le T,\dbP_0\mbox{-a.s.}\nonumber
\end{eqnarray}
We refer to the seminal paper by Pardoux and Peng \cite{PP1} for the
well-posedness of such BSDEs. On the other hand, for any $(t,\omega
)\in\L$, by \cite{PP1} the following BSDE on $[t,T]$ has a unique solution,
%
\begin{eqnarray}
\label{Y0t} Y^{0, t,\omega}_s &=& g^{t,\omega}
\bigl(B^t_\cdot\bigr) + \int_s^T
f^{t,\omega}\bigl(r,B^t_\cdot, Y^{0,t, \omega}_r,
Z^{0,t,\omega}_r\bigr) \,dr \nonumber\\[-8pt]\\[-8pt]
&&{}- \int_s^T
Z^{0,t,\omega}_r\,dB^t_r,\qquad
\dbP_0^t\mbox{-a.s.}\nonumber
\end{eqnarray}
By the Blumenthal 0--1 law, $Y^{0,t,\omega}_t$ is a constant and we
thus define
%
\begin{equation}
\label{u=Y} u^0(t,\omega):= Y^{0,t,\omega}_t.
\end{equation}

\begin{theorem}
\label{thm-existence}
Under Assumption \ref{assum-fg}, $u^0$ is a viscosity solution of PPDE
(\ref{PPDE}) with terminal condition $g$.
\end{theorem}

The proof is reported in Section \ref{sect-existence}.
Similar to the classical theory of viscosity solutions in the
Markovian case, we now establish a comparison result which, in
particular, implies the uniqueness of viscosity solutions. For this
purpose, we need an additional condition:
%
\begin{assum}
\label{assum-hatf}
There exist $\hat f\dvtx  \hat\L\times\dbR\times\dbR^d \to\dbR$ satisfying:

\begin{longlist}[(ii)]
\item[(i)]
$\hat f(t,\omega, y, z) = f(t,\omega, y,z)$ for all $(t,\omega, y,
z) \in\L\times\dbR\times\dbR^d$.

\item[(ii)] $\hat f$ is bounded, $\hat f(\cdot, y,z) \in C^0(\hat\L)$ for
any fixed $(y,z)$ and $\hat f$ is uniformly Lipschitz continuous in $(y,z)$.
\end{longlist}
\end{assum}
%
\begin{rem}
\label{rem-hatf}
In the Markovian case as in Remark \ref{rem-markovian}, we have a
natural extension: $\hat f = f(t, \hat\omega_t, y, z)$ for all $\hat
\omega\in\hat\O$. In this case Assumption \ref{assum-fg} implies
Assumption \ref{assum-hatf}.
\end{rem}
%
\begin{theorem}
\label{thm-comparison}
Let Assumptions \ref{assum-fg} and \ref{assum-hatf} hold.
Let $u^1$ be a viscosity subsolution and $u^2$ a viscosity
supersolution of PPDE (\ref{PPDE}). If $u^1(T,\cdot) \le g\le
u^2(T,\cdot)$, then $u^1 \le u^2$ on $\L$.

Consequently, given the terminal condition $g$, $u^0$ is the unique
viscosity solution of PPDE (\ref{PPDE}).
\end{theorem}

The proof is reported in Section \ref{sect-peron} building on a
partial comparison result derived in Section \ref{sect-partialcomparison}.
%
\begin{rem}
\label{rem-comparison}
For technical reasons, we require a uniformly continuous function
$g$ between $u^1_T$ and $u^2_T$; see Section \ref{sect-peron}.
However, when one of $u^1$ and $u^2$ is in $C^{1,2}_b(\L)$, then we
need neither the presence of such $g$ nor the existence of $\hat f$;
see Lemma \ref{lem-comparison} below.
\end{rem}

\section{Some proofs of the main results}
\label{sect-proof}

In this section we provide some proofs of the main results, and provide
some more results. We leave the most technical part of the proof for
the comparison principle to next section.

\subsection{Properties of classical solutions}
\label{sect-classical}

We first recall from Peng \cite{Peng-mono} that an $\dbF
$-progressively measurable process $Y$ is called an $f$-martingale
(resp., $f$-submartingale, $f$-supermartingale) if, for any $\dbF
$-stopping times $\t_1\ge\t_2$, we have
\[
Y_{\t_1} = (\mbox{resp.,} \le, \ge) \cY_{\t_1}(
\t_2, Y_{\t
_2}),\qquad \dbP_0\mbox{-a.s.},
\]
where $(\cY, \cZ):= ( \cY(\t_2, Y_{\t_2}), \cZ(\t_2, Y_{\t
_2}))$ is the solution to the following BSDE on $[0,\t_2]$:
\[
\cY_t = Y_{\t_2} + \int_t^{\t_2}
f(s, B_\cdot, \cY_s, \cZ_s) \,ds - \int
_t^{\t_2} Z_s \,dB_s,\qquad 0\le t
\le\t_2,\dbP_0\mbox{-a.s.}
\]
Clearly, $Y$ is an $f$-martingale with terminal condition $g(B_\cdot
)$ if and only if it satisfies the BSDE (\ref{BSDE}).

Applying It\^{o}'s formula to Proposition \ref{prop-local-Ito}, we
obviously have the following:

\begin{prop}
\label{prop-classical} Let Assumption \ref{assum-fg} hold and $u\in
C^{1,2}_b(\L)$. Then $u$ is a classical solution (resp., subsolution,
supersolution) of PPDE (\ref{PPDE}) if and only if the process $u$ is
an $f$-martingale (resp., $f$-submartingale, $f$-super\-martingale).

In particular, if $u$ is a classical solution of PPDE (\ref{PPDE})
with terminal condition~$g$, then
%
\begin{equation}
\label{Y=u} Y:= u,\qquad Z:= \partial_{\omega} u
\end{equation}
provides the unique solution of BSDE (\ref{BSDE}).
\end{prop}
\begin{pf} We shall only prove the subsolution case. Let $(Y, Z)$ be
defined by (\ref{Y=u}).

\begin{longlist}[(iii)]
\item[(i)]
Assume $u$ is a classical subsolution. By It\^o's formula,
%
\begin{eqnarray}
\label{uIto} d u_t &=& \bigl(\partial_t u_t
+ \tfrac12 \operatorname{tr} \bigl[\partial^2_{\omega\omega
}
u_t \bigr] \bigr) \,dt + \partial_{\omega} u_t
\,dB_t
\nonumber\\[-8pt]\\[-8pt]
&=& - \bigl(f(t, B_\cdot, u_t, \partial_{\omega}
u_t) + (\cL u)_t \bigr)\,dt + \partial_{\omega}
u_t \,dB_t,\qquad \dbP_0\mbox{-a.s.}\nonumber
\end{eqnarray}
Then for any $\t_1\le\t_2$, $(Y, Z)$ satisfies BSDE
\begin{eqnarray*}
Y_t &=& u_{\t_2} + \int_t^{\t_2}
\bigl(f(s, B_\cdot, Y_s, Z_s) +(\cL u)
(s,B_\cdot) \bigr)\,ds\\
&&{} - \int_t^{\t_2}
Z_s \,dB_s,\qquad 0\le t\le\t_2, \dbP_0
\mbox{-a.s.}
\end{eqnarray*}
Since $\cL u \le0$, by the comparison principle of BSDEs (see \cite
{EPQ}) we obtain $u_{\t_1} = Y_{\t_1} \le\cY_{\t_1}(\t_2, u_{\t
_2})$. That is, $u$ is an $f$-submartingale.

\item[(ii)] Assume $u$ is an $f$-submartingale. For any $0\le t< t+h\le T$,
denote $\d Y_s:= \cY_s(t+h, u_{t+h})-Y_s$, $\d Z_s:= \cZ_s(t+h,
u_{t+h})-Z_s$. By (\ref{uIto}) we have
\begin{eqnarray*}
\d Y_s &=& \int_s^{t+h}\bigl[
\a_r\d Y_r + \langle\b, \d Z\rangle_r - (
\cL u)_r\bigr]\,dr \\
&&{}- \int_s^{t+h} \d
Z_s \,dB_s,\qquad t\le s\le t+h, \dbP_0\mbox{-a.s.},
\end{eqnarray*}
where $|\a|, |\b|\le L_0$. Define
%
\begin{equation}
\label{G} \G_s:= \exp\biggl(\int_t^s
\b_r \,dB_r +\int_t^s
\biggl(\a_r - \frac12 |\b_r|^2 \biggr)\,dr
\biggr),
\end{equation}
and we have
\[
\d Y_t = -\dbE^{\dbP_0}_t \biggl[ \int
_t^{t+h} \G_s(\cL u)_s\,ds
\biggr].
\]
Since $Y=u$ is an $f$-submartingale, we get
\[
0 \le\frac{1}{h} \d Y_t = -\dbE^{\dbP_0}_t
\biggl[ \frac{1}{h}\int_t^{t+h}
\G_s(\cL u)_s\,ds \biggr].
\]
Send $h\to0$, and we obtain
\[
\cL u(t,B_\cdot) \le0, \qquad\dbP_0\mbox{-a.s.}
\]
Note that $\cL u$ is continuous in $\omega$ and obviously any support
of $\dbP_0$ is dense, and we have
\[
\cL u(t,\omega) \le0 \qquad\mbox{for all }\omega\in\O.
\]
That is, $u$ is a classical subsolution of PPDE (\ref{PPDE}).

\item[(iii)] When $u$ is a classical solution similar to (i), we know $Y$ is a
$f$-martingale, and thus (\ref{Y=u}) provides a solution to the BSDE.
Finally, the uniqueness follows from the uniqueness of BSDEs.
\end{longlist}
\upqed\end{pf}
%
\begin{rem}
\label{rem-FK}
This proposition extends the well-known nonlinear Feyn\-man--Kac
formula of Pardoux and Peng \cite{PP2} to non-Markovian case.
\end{rem}

We next prove a simple comparison principle for classical solutions.
%
\begin{lem}
\label{lem-classical-comparison}
Let Assumption \ref{assum-fg} hold true. Let $u^1$ be a classical
subsolution and $u^2$ a classical supersolution of PPDE (\ref{PPDE}).
If $u^1(T,\cdot) \le u^2(T,\cdot)$, then $u^1 \le u^2$ on $\L$.
\end{lem}
\begin{pf} Denote $Y^i:= u^i, Z^i:= \partial_{\omega} u^i$,
$i=1,2$. By
(\ref{uIto}) we have
\[
d Y^i_t = - \bigl[f\bigl(t, B_\cdot,
Y^i, Z^i\bigr)+\bigl(\cL u^i
\bigr)_t \bigr]\,dt + Z^i_t \,dB_t,\qquad 0
\le t\le T, \dbP_0\mbox{-a.s.}
\]
Since $Y^1_T\le Y^2_T$ and $\cL u^1 \le0\le\cL u^2$, by the
comparison principle for BSDEs we obtain $Y^1\le Y^2$. That is, $u^1\le
u^2$, $\dbP_0$-a.s. Since $u^1$ and $u^2$ are continuous, and the
support of $\dbP_0$ is dense in $\O$, we obtain $u^1\le u^2$ on $\L$.
\end{pf}

\subsection{Existence of viscosity solutions}
\label{sect-existence}

We first establish the regularity of $u^0$ as defined in (\ref{u=Y}).
%
\begin{prop}
\label{prop-u0}
Under Assumption \ref{assum-fg}, $u^0$ is uniformly continuous in $\L
$ under $d_\infty$.
\end{prop}
\begin{pf} Since $f$ and $g$ are bounded, clearly $u^0$ is bounded. To
show the uniform continuity, let $(t_i, \omega^i)\in\L$, $i=1,2$,
and assume
without loss of generality that $0\le t_1 \le t_2 \le T$. By taking
conditional expectations $\dbE^{\dbP_0^{t_1}}_{t_2}$, one can easily
see that $Y^{0, t_1, \omega^1}$ can be viewed as the solution to the
following BSDE on $[t_2, T]$: for $\dbP^{t_1}_0$-a.s. $B^{t_1}$,
\begin{eqnarray*}
Y^{0, t_1,\omega^1}_s &=& g^{t_2, \omega^1\otimes_{t_1}
B^{t_1}}\bigl(B^{t_2}_\cdot
\bigr) + \int_s^T f^{t_2,\omega^1\otimes_{t_1}
B^{t_1}}
\bigl(r,B^{t_2}_\cdot, Y^{0,t_1,\omega^1}_r,
Z^{0,t_1,\omega
^1}_r\bigr) \,dr
\\
&&{}- \int_s^T Z^{0,t_1,\omega^1}_r\,dB^{t_2}_r,\qquad t_2
\le s\le T,\dbP_0^{t_2}\mbox{-a.s.}
\end{eqnarray*}
Denote
\[
\d\omega:=\omega^1-\omega^2,\qquad \d Y:= Y^{0, t_1,\omega^1}-
Y^{0, t_2,\omega^2},\qquad \d Z:= Z^{0,t_1,\omega^1}-Z^{0,t_2,\omega^2}.
\]
Then
\begin{eqnarray*}
\d Y_s &=& \d Y_T + \int_s^T
\bigl(\gamma_r + \a_r \d Y_r + \langle\b,
\d Z\rangle_r \bigr) \,dr \\
&&{}- \int_s^T
\d Z_r \,dB^{t_2}_r,\qquad t_2\le s\le T,
\dbP_0^{t_2}\mbox{-a.s.},
\end{eqnarray*}
where
\[
|\a|\le L_0,\qquad \b\in\cU^{L_0}_t
\]
and
\[
\gamma_r:= \bigl(f^{t_2,\omega^1\otimes_{t_1}
B^{t_1}}-f^{t_2,\omega^2} \bigr)
\bigl(r,B^{t_2}_\cdot,Y^{0, t_1,\omega^1},Z^{0, t_1,\omega
^1}_r
\bigr).
\]
Define $\G$ as in (\ref{G}) with initial time $t_2$,
then
\[
\d Y_{t_2} = \G_T \d Y_T + \int
_{t_2}^T \G_r \gamma_r
\,dr - \int_{t_2}^T \G_r[\d
Z_r + \d Y_r \b_r] \,dB^{t_2}_r,\qquad
\dbP_0^{t_2}\mbox{-a.s.}
\]
Let $\rho$ denote the modulus of continuity function of $f$ and $g$
with respect to $\omega$. Note that
\begin{eqnarray*}
|\d Y_T| &=& \bigl| g\bigl(\omega^1\otimes
_{t_1} B^{t_1}\otimes
_{t_2} B^{t_2}\bigr) - g\bigl(\omega^2
\otimes_{t_2} B^{t_2}\bigr) \bigr|
\\
&\le& \rho\bigl( \bigl\|\omega^1\otimes
_{t_1} B^{t_1} - \omega^2
\bigr\|_{t_2} \bigr) \le\rho\bigl(d_\infty\bigl(
\bigl(t_1, \omega^1\bigr), \bigl(t_2,
\omega^2\bigr) \bigr) + \bigl\|B^{t_1}\bigr\|^{t_1}_{t_2}
\bigr).
\end{eqnarray*}
Similarly,
\[
|\gamma_s|\le\rho\bigl(d_\infty\bigl(\bigl(t_1,
\omega^1\bigr), \bigl(t_2,\omega^2\bigr)
\bigr) + \bigl\|B^{t_1}\bigr\|^{t_1}_{t_2} \bigr).
\]
Then
\[
|\d Y_{t_2}| = \biggl|\dbE^{\dbP_0^{t_2}} \biggl[ \G_T \d
Y_T + \int_{t_2}^T \G_s
\gamma_s \,ds \biggr] \biggr| \le C\rho\bigl(d_\infty\bigl(
\bigl(t_1, \omega^1\bigr), \bigl(t_2,
\omega^2\bigr) \bigr) + \bigl\|B^{t_1}\bigr\|^{t_1}_{t_2}
\bigr).
\]
Thus, noting that $f$ is bounded,
%
\begin{eqnarray}
\label{Du}
&&\bigl|u^0_{t_1}\bigl(\omega^1\bigr) -
u^0_{t_2}\bigl(\omega^2\bigr)\bigr| \nonumber\\
&&\qquad=
\bigl|Y^{0,t_1, \omega^1}_{t_1} - Y^{0,t_2, \omega^2}_{t_2}\bigr|
\nonumber
\\
&&\qquad= \biggl|\dbE^{\dbP_0^{t_1}} \biggl[Y^{0,t_1, \omega^1}_{t_2} + \int
_{t_1}^{t_2} f^{t_1,\omega^1}\bigl(s,B^{t_1}_\cdot,
Y^{0,t_1,\omega
^1}_s, Z^{0,t_1,\omega^1}_s\bigr) \,ds-
Y^{0,t_2, \omega^2}_{t_2} \biggr]\biggr|
\\
&&\qquad\le C[t_2-t_1] + \dbE^{\dbP_0^{t_1}} \bigl[|\d
Y_{t_2}| \bigr]
\nonumber
\\
&&\qquad\le C[t_2-t_1]+ C\dbE^{\dbP_0^{t_1}} \bigl[\rho
\bigl(d_\infty\bigl(\bigl(t_1, \omega^1\bigr),
\bigl(t_2,\omega^2\bigr) \bigr) + \bigl\|B^{t_1}
\bigr\|^{t_1}_{t_2} \bigr) \bigr].\nonumber
\end{eqnarray}
For any $\e>0$, there exists $h>0$ such that $\rho(h) \le\frac{\e
}{2C}$ for the above $C$. Since $f, g$ are bounded, we may assume $\rho
$ is also bounded and denote by $\|\rho\|_\infty$ its bound. Now for
$ d_\infty((t_1, \omega^1), (t_2,\omega^2) )\le\frac
{h}{2}$, we obtain
\begin{eqnarray*}
&&
\bigl|u^0_{t_1}\bigl(\omega^1\bigr) -
u^0_{t_2}\bigl(\omega^2\bigr)\bigr| \\
&&\qquad\le
Cd_\infty\bigl(\bigl(t_1, \omega^1\bigr),
\bigl(t_2,\omega^2\bigr) \bigr) + C\rho(h) + C\|\rho
\|_\infty\dbP_0^{t_1} \biggl[\bigl\|B^{t_1}
\bigr\|^{t_1}_{t_2} > \frac
{h}{2} \biggr]
\\
&&\qquad\le \frac{\e}{2} + Cd_\infty\bigl(\bigl(t_1,
\omega^1\bigr), \bigl(t_2,\omega^2\bigr)
\bigr) + 4C\|\rho\|_\infty h^{-2} \dbE^{\dbP_0^{t_1}} \bigl[
\bigl(\bigl\|B^{t_1}\bigr\|^{t_1}_{t_2}\bigr)^2
\bigr]
\\
&&\qquad= \frac{\e}{2} + Cd_\infty\bigl(\bigl(t_1,
\omega^1\bigr), \bigl(t_2,\omega^2\bigr)
\bigr) + 4C\|\rho\|_\infty h^{-2}(t_2-t_1)
\\
&&\qquad\le \frac{\e}{2} + C \bigl(1+4\|\rho\|_\infty h^{-2}
\bigr) d_\infty\bigl(\bigl(t_1, \omega^1\bigr),
\bigl(t_2,\omega^2\bigr) \bigr).
\end{eqnarray*}
By choosing $d_\infty((t_1, \omega^1), (t_2,\omega^2) )$
small enough, we see that $|u^0_{t_1}(\omega^1) -\break  u^0_{t_2}(\omega
^2)|\le\e$. This completes the proof.
\end{pf}

However, in general one cannot expect $u^0$ to be a classical solution
to PPDE (\ref{PPDE}). We refer to Peng and Wang \cite{PW} for some
sufficient conditions, in a slightly different setting.
\begin{pf*}{Proof of Theorem \ref{thm-existence}}
We just show that $u^0$ is a viscosity subsolution. We prove by
contradiction. Assume $u^0$ is not a viscosity subsolution. Then, for
all $L>0$, $u^0$ is not an $L$-viscosity subsolution. For the purpose
of this proof, it is sufficient to consider an arbitrary $L\ge L_0$,
the Lipschitz constant of $f$ introduced in Assumption \ref{assum-fg}(i). Then, there exist
\[
(t,\omega)\in[0,T)\times\O\quad\mbox{and}\quad \f\in
\underline
\cA^{L} u^0 (t,\omega) \qquad\mbox{such that } c:= \bigl(
\cL^{t,\omega}\f\bigr) (t,0)>0.
\]
Denote, for $s\in[t,T]$,
\begin{eqnarray*}
\tilde Y_s &:=& \f\bigl(s,B^t\bigr),\qquad \tilde
Z_s:= \partial_{\omega}\f\bigl(s,B^t\bigr),\\
\d Y_s&:=&\tilde Y_s - Y^{0,t,\omega}_s,\qquad \d
Z_s:=\tilde Z_s- Z^{0,t,\omega}_s.
\end{eqnarray*}
Applying It\^{o}'s formula, we have
\begin{eqnarray*}
d (\d Y_s) &=& - \bigl[\bigl(\cL^{t,\omega} \f\bigr)
\bigl(s, B^t_\cdot\bigr) + f^{t,\omega}\bigl(s,
B^t_\cdot, \tilde Y_s, \tilde Z_s
\bigr) - f^{t,\omega
}\bigl(s, B^t_\cdot,
Y^{0,t,\omega}_s, Z^{0,t,\omega}_s\bigr) \bigr]\,ds\\
&&{}+ \d
Z_s \,dB^t_s
\\
&=& - \bigl[\bigl(\cL^{t,\omega} \f\bigr) \bigl(s, B^t_\cdot
\bigr) + \a_s \d Y_s + \langle\b, \d Z
\rangle_s \bigr]\,ds + \d Z_s \,dB^t_s,\qquad
\dbP_0^t\mbox{-a.s.},
\end{eqnarray*}
where $|\a|\le L_0$ and $\b\in\cU^{L_0}_t\subset\cU^L_t$.
Observing that $\d Y_t=0$, we define
\[
\t_0:= T\wedge\inf\biggl\{s > t\dvtx  \bigl(
\cL^{t,\omega} \f\bigr) \bigl(s, B^t_\cdot\bigr) -
L_0 |\d Y_s|\le\frac{c}{2} \biggr\}.
\]
Then, by Proposition \ref{prop-u0} and Example \ref{eg-tau}, $\t_0\in\cT
^t_+$ and
%
\begin{equation}
\label{Lphi<0} \bigl(\cL^{t,\omega} \f\bigr) \bigl(s, B^t_\cdot
\bigr) + \a_s \d Y_s\ge\frac
{c}{2}\qquad\mbox{for all } s
\in[t,\t_0].
\end{equation}
Now for any $\t\in\cT^t$ such that $\t\le\t_0$, we have
\begin{eqnarray*}
0 &=& \d Y_t = \d Y_{\t} + \int_t^{\t}
\bigl[ \bigl(\cL^{t,\omega} \f\bigr) \bigl(s, B^t_\cdot
\bigr) + \a_s \d Y_s + \langle\b, \d Z
\rangle_s \bigr]\,ds - \int_t^{\t} \d
Z_s \,dB^t_s
\\
&\ge& \f\bigl(\t, B^t\bigr) - u^{0,t,\omega}\bigl(\t,
B^t\bigr) + \frac{c}{2}(\t-t) - \int_t^{\t}
\d Z_s \bigl(dB^t_s-\beta_s\,ds
\bigr).
\end{eqnarray*}
Then\vspace*{1pt} $\underline\cE^L_t [(\f- u^{0,t,\omega})(\t, B^t)
]\le
\dbE^{\dbP^\beta}_t [(\f- u^{0,t,\omega})(\t, B^t) ]\le
0$. This contradicts with $\f\in\underline\cA^{L} u^0(t,\omega)$.
\end{pf*}

Following similar arguments, one can easily prove the following:
%
\begin{prop}
\label{prop-classical-viscosity}
Under Assumption \ref{assum-fg}, a bounded classical subsolution
(resp., supersolution) of the PPDE (\ref{PPDE}) must be a viscosity
subsolution (resp., supersolution).
\end{prop}

\subsection{Stability of viscosity solutions}
\label{sect-stability}

\mbox{}

\begin{pf*}{Proof of Theorem \ref{thm-stability}} We shall prove only
the viscosity subsolution property by contradiction. By Remark \ref
{rem-strict}, without loss of generality we assume there exists $\f\in
{\underline\cA'}^{L}u(0,{\mathbf0})$ such that $c:= \cL\f(0,{\mathbf0})
> 0$, where ${\underline\cA'}^{L}u(0,{\mathbf0})$ is defined in~(\ref{cA1}).

Denote
%
\begin{eqnarray}
\label{tau0}
X^0&:=& \f- u,\qquad X^\e:= \f- u^{\e}
\quad\mbox{and}\nonumber\\[-8pt]\\[-8pt]
\t_0&:=& \inf\biggl\{ t>0\dvtx  \cL\f(t, B) \le\frac{c}{2}
\biggr\}\wedge T.\nonumber
\end{eqnarray}
Since $f\in C^0(\Lambda)$, it follows from Example \ref{eg-tau} that
$\t_0\in\cT^0_+$. By (\ref{cA1}), there exists $\t_1\in\cT^0_+$
such that $\t_1 \le\t_0$ and
\[
\underline\cE^L_0\bigl(\t_1,
X^0_{\t_1}\bigr) > 0 = X^0_0.
\]
Since $u^\e$ converges toward $u$ uniformly, we have
%
\begin{equation}
\label{f-u>0} \underline\cE^L_0\bigl(\t_1,
X^\e_{\t_1}\bigr) > X^\e_0 \qquad\mbox{for
sufficiently small } \e>0.
\end{equation}
Consider the optimal stopping problem, under nonlinear expectation,
together with the corresponding optimal stopping rule,
%
\begin{equation}
Y_t := Y_t^\e := \inf_{\tau\in\cT^t} \underline{\cE}^L_t
\bigl[X^\e_{\tau\wedge\tau_1} \bigr] \quad\mbox{and}\quad \t^*_0 :=
\inf\bigl\{t\ge0\dvtx  Y_t = X^\e_t
\bigr\};
\end{equation}
see Remark \ref{rem-cEreflected}. We claim that
%
\begin{equation}
\dbP_0 \bigl[\t^*_0<\t_1 \bigr] > 0,
\end{equation}
because otherwise $X_0^\e\ge Y_0 = \underline{\cE}^L_0 [X^\e_{\t_1} ]$,
contradicting (\ref{f-u>0}).\vspace*{1pt}

Since $X^\e$ and $Y$ are continuous, $\dbP_0$-a.s. there exists $E\subset
\{\t^*_0<\t_1\}$ such that $\dbP_0(E) = \dbP_0(\t^*_0<\t_1)>0$,
and for any $\omega\in E$, denoting $t:=\t^*_0(\omega)$ we have
$X_t^\e (\omega) = Y_t(\omega)$. Notice that $\t^{t,\omega}_1 \in
\cT^t_+$. By standard arguments using the regular conditional
probability distributions (see, e.g., \cite{SV} or \cite
{STZ-duality}), it follows from the definition of $\tau^*_0$ together
with the $\cE^L$-submartingale property of $Y$ that
\[
X^\e_t(\omega) = Y_t(\omega) =
Y^{t,\omega}_t(\omega) \le\underline{\cE}^L_t
\bigl[Y^{t,\omega}_\t\bigr] \le\underline{\cE}^L_t
\bigl[X^{\e,t,\omega}_\t\bigr] \qquad\mbox{for all } \tau\in
\cT^t, \t\le\t^{t,\omega}_1.
\]
This implies that
\begin{eqnarray*}
0 &\le&\underline\cE^{L}_t \bigl[X^{\e,t,\omega}_{\t}
- X^\e_t(\omega) \bigr] \\
&=& \underline
\cE^{L}_t \bigl[\f^{t,\omega}_\t-\f(t,
\omega) + u^\e(t,\omega)-u^{\e, t,\omega}_\t\bigr]
\qquad\mbox{for all } \tau\in\cT, \t\le\t^{t,\omega}_1.
\end{eqnarray*}
Define
\[
\f^\e_s:= \f^{t,\omega}_s - \f(t,\omega)
+ u^\e(t,\omega).
\]
Then we have $\f^\e\in\underline\cA^Lu^\e(t,\omega)$. Since
$u^\e$ is a viscosity $L$-subsolution of PPDE (\ref{PPDE}) with
coefficients $f^\e$, we have
\begin{eqnarray*}
0 &\ge& -\partial_t \f^\e(t,0) - \frac12
\operatorname{tr} \bigl[\partial^2_{\omega\omega}\f^\e
\bigr](t,0) - f^\e\bigl(t,\omega, \f^\e(t,0),
\partial_{\omega}\f^\e(t,0) \bigr)
\\
&=& -\partial_t \f(t,\omega) - \frac12\operatorname{tr} \bigl[
\partial^2_{\omega\omega}\f\bigr](t,\omega) - f^\e
\bigl(t,\omega, u^\e(t,\omega), \partial_{\omega}\f(t,\omega)
\bigr)
\\
&=& (\cL\f) (t,\omega) + f \bigl(t,\omega, u(t,\omega), \partial_{\omega}
\f(t,\omega) \bigr) - f^\e\bigl(t,\omega, u^\e(t,
\omega), \partial_{\omega}\f(t,\omega) \bigr)
\\
&\ge& \frac{c}{2} + f \bigl(t,\omega, u(t,\omega), \partial_{\omega}
\f(t,\omega) \bigr) - f^\e\bigl(t,\omega, u^\e(t,
\omega), \partial_{\omega}\f(t,\omega) \bigr),
\end{eqnarray*}
thanks to (\ref{tau0}).
Send $\e\to0$, we obtain $0\ge\frac{c}{2}$, a contradiction.
\end{pf*}
%
\begin{rem}
\label{rem-stability}
(i) We need the same $L$ in the proof of Theorem \ref
{thm-stability}. If $u^\e$ is only a viscosity subsolution of PPDE
(\ref{PPDE}) with coefficient $f^\e$, but with possibly different
$L_\e$, we are not able to show that $u$ is a viscosity subsolution of
PPDE (\ref{PPDE}) with coefficients (\ref{PPDE}).

(ii) However, if $u^\e$ is a viscosity solution of PPDE (\ref{PPDE})
with coefficient $f^\e$, by Theorems \ref{thm-existence} and \ref
{thm-comparison}, it follows immediately from the stability of BSDEs
that $u$ is the unique viscosity solution of PPDE (\ref{PPDE}) with
coefficient $f$.
\end{rem}
%

\subsection{Partial comparison principle}
\label{sect-partialcomparison}

The following partial comparison principle, which improves Lemma \ref
{lem-classical-comparison}, is crucial for this paper. The main
argument is very much similar to that of Theorem \ref{thm-stability}.
%
\begin{lem}
\label{lem-comparison}
Let Assumption \ref{assum-fg} hold true. Let $u^1$ be a viscosity
subsolution and $u^2$ a viscosity supersolution of PPDE (\ref{PPDE}).
If $u^1(T,\cdot) \le u^2(T,\cdot)$ and one of $u^1$ and $u^2$ is
in $C^{1,2}_b(\L)$, then $u^1 \le u^2$ on $\L$.
\end{lem}
\begin{pf} First, by Remark \ref{rem-change}(i), by otherwise
changing the variable we may assume without loss of generality that
%
\begin{equation}
\label{fmonotone}
\mbox{$f$ is strictly decreasing in $y$.}
\end{equation}
We assume $u^2 \in C^{1,2}_b(\L)$ and $u^1$ is a viscosity
$L$-subsolution for some $L\ge0$. We shall prove by contradiction.
Without loss of generality, we assume
\[
-c:= u^2_0- u^1_0 < 0.
\]
For future purposes, we shall obtain the contradiction under the
following slightly weaker assumptions:
%
\begin{eqnarray}
\label{u2-weak} u^2 &\in&\bar C^{1,2}_{\dbP_0}(\L)\qquad
\mbox{bounded}\quad \mbox{and}\nonumber\\[-8pt]\\[-8pt]
\bigl(\cL u^2\bigr) &\ge&0,\qquad
u^2(T,\cdot) \ge u^1(T,\cdot) ,\qquad\dbP_0
\mbox{-a.s.}\nonumber
\end{eqnarray}
Denote
\[
X:= u^2-u^1 \quad\mbox{and}\quad \t_0:= \inf\{ t>0\dvtx
X_t \ge0 \} \wedge T.
\]
Note that $X_0 = -c <0$, $X_T\ge0$, and $X$ is continuous, $\dbP
_0$-a.s. Then
%
\begin{equation}
\label{X<0}\quad \t_0>0,\qquad X_t < 0,\qquad t\in[0, \t_0),
\quad\mbox{and}\quad X_{\t
_0} = 0,\qquad \dbP_0\mbox{-a.s.}
\end{equation}
Similar to Remark \ref{rem-cEreflected}, define the process $Y$ by
the optimal stopping problem under nonlinear expectation,
\[
Y_t := \inf_{\tau\in\cT^t}\underline{\cE}^L_t
[X_{\tau\wedge\tau
_0} ],\qquad t\in[0,\tau_0],
\]
together with the corresponding optimal stopping rule,
\[
\tau^*_0 := \inf\{t\ge0\dvtx Y_t=X_t\}.
\]
Then $\tau^*_0\le\tau_0$, and we claim that
%
\begin{equation}
\dbP_0\bigl[\tau^*_0<\tau_0\bigr] > 0,
\end{equation}
because otherwise $X_0\ge Y_0=\underline{\cE}^L_0 [X_{\tau
_0} ]$, contradicting (\ref{X<0}).

As in the proof of Theorem \ref{thm-stability}, there exists $E\subset
\{\t^*_0<\t_0\}$ such that $\dbP_0(E) = \dbP_0 [\t^*_0<\t_0 ]>0$, and
for any $\omega\in E$, by denoting $t:=\t^*_0(\omega)$ we have $\tau
^{t,\omega}_0\in\cT^{t}_+$ and
\[
X_t(\omega) = Y_t(\omega) = \inf_{\tau\in\cT^t}
\underline{\cE}^L_t [X_{\tau\wedge\t
^{t,\omega}_0}],\qquad
\dbP^t_0\mbox{-a.s.}
\]
Let $\{\t_i, i\ge0\}$ be the sequence of stopping times in Definition
\ref{defn-barCP} corresponding to $u^2$. Then $\dbP_0 [\{\t^*_0
< \t_i\}\cap E ]>0$ for $i$ large enough, and thus there exists
$\omega\in E$ such that $t:= \t^*_0(\omega) < \t_i(\omega)$.
Without loss of generality, we assume $\t_{i-1}(\omega) \le t$. It
is clear that $(\t_0\wedge\t_{i})^{t,\omega} \in\cT^t_+$ and
$(u^2)^{t,\omega} \in C^{1,2}_b(\bar\L^t((\t_0\wedge\t_{i})^{t,\omega
}))$. In particular, there exists $\tilde u^2\in
C^{1,2}_b(\L^t)$ such that $(u^2)^{t,\omega} = \tilde u^2$ on $(\t
_0\wedge\t_{i})^{t,\omega}$.

Now for any $\t\in\cT^t_+$ such that $\t\le(\t_0\wedge\t_{i})^{t,\omega
}$, it follows from Remark \ref{rem-cEreflected} that
\[
X_t(\omega) = Y_t(\omega) =Y^{t,\omega}_t
\le\underline\cE^{L}_t \bigl[Y^{t,\omega}_{\t}
\bigr] \le\underline\cE^{L}_t \bigl[X^{t,\omega}_{\t}
\bigr].
\]
Thus
\[
0\le\underline\cE^{L}_t \bigl[\bigl(\tilde
u^2\bigr)^{t,\omega}_{\t} - \bigl(u^1
\bigr)^{t,\omega}_{\t} - X_t(\omega) \bigr].
\]
Denote $\f_s:= (\tilde u^2)^{t,\omega}_s - X_t(\omega)$, $s\in
[t,T]$. Then $\f\in\underline\cA^Lu^1(t,\omega)$. Since $u^1$ is
a viscosity $L$-subsolution, and $u^2$ is a classical supersolution, we have
\begin{eqnarray*}
0 &\ge& (\cL\f) (t,0) = - \partial_t \tilde u^2(t,0) -
\tfrac12 \operatorname{tr} \bigl[\partial^2_{\omega\omega} \tilde
u^2(t,0) \bigr] - f \bigl(t,\omega, u^1(t,\omega),
\partial_{\omega}\tilde u^2(t,0) \bigr)
\\
&=& - \partial_t u^2(t,\omega) - \tfrac12
\operatorname{tr} \bigl[\partial^2_{\omega\omega
}u^2(t,
\omega) \bigr] - f \bigl(t,\omega, u^1(t,\omega),
\partial_{\omega}u^2(t,\omega) \bigr)
\\
&=& \bigl(\cL u^2\bigr) (t,\omega) + f \bigl(t,\omega,
u^2(t,\omega), \partial_{\omega}u^2(t,\omega)
\bigr) - f \bigl(t,\omega, u^1(t,\omega), \partial_{\omega}u^2(t,
\omega) \bigr)
\\
&\ge& f \bigl(t,\omega, u^2(t,\omega), \partial_{\omega}u^2(t,
\omega) \bigr) - f \bigl(t,\omega, u^1(t,\omega),
\partial_{\omega}u^2(t,\omega) \bigr).
\end{eqnarray*}
By (\ref{X<0}), $u^2(t,\omega) < u^1(t,\omega)$. Then the above
inequality contradicts with (\ref{fmonotone}).\vadjust{\goodbreak}
\end{pf}

\section{A variation of Perron's approach}
\label{sect-peron}

To prove Theorem \ref{thm-comparison}, we define
%
\begin{eqnarray}
\label{baru} \overline u(t,\omega)&:=& \inf\bigl\{\f(t,{\mathbf0})\dvtx  \f\in
\overline
\cD(t,\omega) \bigr\},\nonumber\\[-8pt]\\[-8pt]
\underline u(t,\omega)&:=& \sup\bigl\{\f(t,{\mathbf
0})\dvtx  \f\in
\underline\cD(t,\omega) \bigr\},\nonumber
\end{eqnarray}
where, in light of (\ref{u2-weak}),
%
\begin{eqnarray}
\overline\cD(t,\omega)&:=& \bigl\{ \f\in\bar C^{1,2}_{\dbP^t_0}
\bigl(\L^t\bigr) \mbox{ bounded: }(\cL\f)^{t,\omega}_s
\ge0,\nonumber\\
&&\hspace*{24pt}s\in[t,T]\mbox{ and } \f_T \ge g^{t,\omega},
\dbP_0^t\mbox{-a.s.} \bigr\};
\nonumber\\[-8pt]\\[-8pt]
\underline\cD(t,\omega)&:=& \bigl\{ \f\in\bar C^{1,2}_{\dbP
^t_0}
\bigl(\L^t\bigr)\mbox{ bounded: }(\cL\f)^{t,\omega}_s
\le0,\nonumber\\
&&\hspace*{24pt}s\in[t,T]\mbox{ and }\f_T \le g^{t,\omega},
\dbP_0^t\mbox{-a.s.} \bigr\}.\nonumber
\end{eqnarray}
By Lemma \ref{lem-comparison}, in particular by its proof under the
weaker condition (\ref{u2-weak}), it is clear that
%
\begin{equation}
\label{underlineu<overlineu} \underline u \le u^0 \le\overline u.
\end{equation}
The following result is important for our proof of Theorem \ref
{thm-comparison}.
%
\begin{theorem}
\label{thm-peron}
Let Assumptions \ref{assum-fg} and \ref{assum-hatf} hold true. Then
%
\begin{equation}
\label{peron} \underline u = \overline u.
\end{equation}
\end{theorem}
\begin{pf*}{Proof of Theorem \ref{thm-comparison}} By Lemma \ref
{lem-comparison}, in particular by its proof under the weaker condition
(\ref{u2-weak}), we have $u^1 \le\overline u$ and $\underline u \le
u^2$. Then Theorem \ref{thm-peron} implies that $u^1\le u^2$.

This clearly leads to the uniqueness of viscosity solution, and
therefore, by Theorem \ref{thm-existence} $u^0$ is the unique
viscosity solution of PPDE (\ref{PPDE}) with terminal condition~$g$.
\end{pf*}
%
\begin{rem}
\label{rem-peron}
In standard Perron's method, one shows that $\underline u$ (resp.,~$\overline u$)
is a viscosity super-solution (resp., viscosity
sub-solution) of the PDE. Assuming that the comparison principle for
viscosity solutions holds true, then (\ref{peron}) holds.

In our situation, we shall instead prove (\ref{peron}) directly first,
which in turn is used to prove the comparison principle for viscosity
solutions. Roughly speaking, the comparison principle for viscosity
solutions is more or less equivalent to the partial comparison
principle Lemma \ref{lem-comparison} and the equality (\ref{peron}).
To our best knowledge, such an approach is novel in the literature.
\end{rem}

We decompose the proof of Theorem \ref{thm-peron} into several lemmas.
First, let $t<T$ and $\th\in(C_b^0(\L^t) )^d$ satisfy
%
\begin{eqnarray}
\label{Zproperty}
\begin{tabular}{p{320pt}}
there exists $\hat\th\in(C_b^0
(\hat\L^t) )^d$ such that $\th= \hat\th$ in
$\L$
and $\hat\th$ is uniformly continuous in $\hat\omega$ under the uniform
norm $\|\cdot\|^t_T$.
\end{tabular}\hspace*{-32pt}
\end{eqnarray}
Define
%
\begin{equation}
\label{v} Z_s = z + \int_t^s
\th_r \,dr,\qquad v_s:= \int_t^s
Z_s \,dB_s^t,\qquad t\le s\le T,\dbP_0^t
\mbox{-a.s.}\hspace*{-28pt}
\end{equation}
By It\^{o}'s formula, we have
\[
v_s = Z_s B^t_s - \int
_t^s \th_r B^t_r
\,dr.
\]
Denote
%
\begin{eqnarray}
\label{hatv} \hat Z_s(\hat\omega)&:=& z + \int_t^s
\hat\th_r(\hat\omega) \,dr,\nonumber\\[-9pt]\\[-9pt]
\hat v(s,\hat\omega)&:=& \hat
Z_s(\hat\omega) \hat\omega_s - \int_t^s
\hat\th_r(\hat\omega) \hat\omega_r \,dr,\qquad \hat\omega\in
\hat\O^t.\nonumber
\end{eqnarray}
Now for any $\omega\in\O$ and $x\in\dbR$, let $\hat
u^{t,\omega}$ denote the unique solution to the following ODE (with
random coefficients) on $[t, T]$:
%
\begin{eqnarray}
\label{Xt} \hat u^{t,\omega}(s,\hat\omega)&:=& x - \int_t^s
\hat f^{t,\omega
}\bigl(r,\hat\omega, \hat u^{t,\omega}(r,\hat\omega), \hat
Z_r(\hat\omega)\bigr) \,dr \nonumber\\[-8pt]\\[-8pt]
&&{}+ \hat v(s,\hat\omega),\qquad t\le s\le T, \hat
\omega\in\hat\O^t,\nonumber
\end{eqnarray}
and define
%
\begin{equation}
\label{uto} u^{t,\omega}(s,\tilde\omega):= \hat u^{t,\omega}(s,\tilde
\omega) \qquad\mbox{for } (s,\tilde\omega)\in\L^t.
\end{equation}
We then have the following:
%
\begin{lem}
\label{lem-Xsmooth}
Let Assumptions \ref{assum-fg} and (\ref{Zproperty}) hold true. Then
for each \mbox{$(t,\omega)\in\L$}, the above $u^{t,\omega}\in
C^{1,2}_b(\L^t)$ and $\cL^{t,\omega} u^{t,\omega} =0$.
\end{lem}
\begin{pf} We first show that $\hat u^{t,\omega} \in C^{1,2}_b(\hat
\L^t)$, which implies that $u^{t,\omega}\in C^{1,2}_b(\L^t)$. For
$t\le s_1<s_2\le T$ and $\hat\omega^1, \hat\omega^2\in\hat\O^t$, we have
\begin{eqnarray*}
\bigl|\hat Z_{s_1}\bigl(\hat\omega^1\bigr)-\hat
Z_{s_2}\bigl(\hat\omega^2\bigr)\bigr| &\le& \int
_{s_1}^{s_2} \bigl|\hat\th_r\bigl(\hat
\omega^1\bigr)\bigr|\,dr + \int_t^{s_1}\bigl|\hat
\th_r\bigl(\hat\omega^1\bigr)-\hat\th_r
\bigl(\hat\omega^2\bigr)\bigr|\,dr
\\[-2pt]
&\le& C[s_2-s_1] + \int_t^{s_1}\bigl|
\hat\th_r\bigl(\hat\omega^1\bigr)-\hat\th_r
\bigl(\hat\omega^2\bigr)\bigr|\,dr.
\end{eqnarray*}
Note that $ d^t_\infty((r,\hat\omega^1), (r,\hat
\omega^2) )\le d^t_\infty((s_1,\hat\omega^1), (s_2,\hat\omega^2) )$ for
$t\le r\le s_1$. Then one can easily see that $\hat Z\in
C_b^0(\hat\L^t)$. Similarly one can show that $\hat v, \hat
u^{t,\omega}\in C^0(\hat\L^t)$.

Next, one can easily check that, for all $\hat\omega\in\hat\O^t$,
\begin{eqnarray*}
\partial_t \hat Z_s (\hat\omega)&=& \hat
\th_s(\hat\omega), \qquad\partial_{\omega} \hat Z_s(\hat
\omega) ={\mathbf0};
\\[-2pt]
\partial_t \hat v (s,\hat\omega)&=& \hat\th_s(\hat
\omega) \hat\omega_s- \hat\th_s(\hat\omega)\hat
\omega_s =0,\\ \partial_{\omega}\hat v (s,\hat\omega)&=& \hat
Z_s(\hat\omega),\qquad
\partial^2_{\omega\omega}\hat v(s,\hat
\omega) = {\mathbf0};\\
\qquad
\partial_t \hat u^{t,\omega}(s,\hat\omega) &=& -\hat
f^{t,\omega
}\bigl(s,\hat\omega, \hat u^{t,\omega}(s,\hat\omega), \hat
Z_s (\hat\omega)\bigr), \\
\partial_{\omega}\hat
u^{t,\omega}(s,\hat\omega) &=& \hat Z_s(\hat\omega),\qquad
\partial^2_{\omega\omega}\hat u^{t,\omega
}(s,\hat\omega) = {
\mathbf0}.
\end{eqnarray*}
Since $\hat\th$ and $\hat f$ are bounded, it is straightforward to
see that $\hat u^{t,\omega} \in C^{1,2}_b(\hat\L^t)$.

Finally, from the above derivatives we see immediately that $\cL
^{t,\omega} u^{t,\omega}=0$.
\end{pf}

Our next two lemmas rely heavily on the remarkable result Bank and Baum
\cite{BB}, which is extended to BSDE case in \cite{STZ-duality}.
%
\begin{lem}
\label{lem-BB1}
Let\vspace*{2pt} Assumption \ref{assum-fg} hold true. Let $\t\in\cT$, $Z$ be
$\dbF$-progressively measurable such that $\dbE^{\dbP_0}[\int_\t^T
|Z_s|^2\,ds] <\infty$, and $X_\t, \tilde X_\t\in L^2(\cF_\t, \dbP_0)$.
For any $\e>0$, there exists $\dbF$-progressively measurable
process $Z^\e$ such that:

\begin{longlist}
\item For the Lipschitz constant $L_0$ in Assumption
\ref{assum-fg}\textup{(ii)}, it holds that
%
\begin{equation}
\label{BBest1} \dbP_0 \Bigl[\sup_{\t\le t\le T}
e^{-L_0t}\bigl|X^{\e}_t - X_t\bigr| \ge\e+
e^{-L_0\t}|\tilde X_\t- X_{\t}| \Bigr] \le\e,
\end{equation}
where $X, X^\e$ are the solutions to the following ODEs with random
coefficients,
%
\begin{eqnarray}
\label{Xe} X_t &=& X_\t- \int_\t^t
f(s,B, X_s, Z_s) \,ds + \int_\t^t
Z_s \,dB_s,
\nonumber
\\
X^\e_t &=& \tilde X_\t- \int
_\t^t f\bigl(s,B, X^\e_s,
Z^\e_s\bigr) \,ds \\
&&{}+ \int_\t^t
Z^\e_s \,dB_s, \qquad\t\le t\le T,\dbP_0
\mbox{-a.s.};\nonumber
\end{eqnarray}

\item  $\th^\e_t:= \frac{d}{dt} Z^\e_t$ exists for $t\in[\t, T)$,
where $\th^\e_\t$ is understood as the right derivative, and for
each $\omega\in\O$, $(\th^\e)^{\t(\omega), \omega}$
satisfies (\ref{Zproperty}) with $t:= \t(\omega)$.
\end{longlist}
\end{lem}
\begin{pf} First, let $h:= h_\e>0$ be a small number which will be
specified later. By standard arguments there exists a time partition
$0=t_0<\cdots<t_n=T$ and a smooth function $\psi\dvtx  [0,T]\times\dbR
^{n\times d}\to\dbR^d$ such that $\psi$ and its derivatives are
bounded and
%
\begin{eqnarray}
\label{tildeZ} &&\dbE^{\dbP_0} \biggl[\int_\t^T
|\tilde Z_t - Z_t|^2\,dt \biggr] <h\hspace*{70pt}
\nonumber\\[-8pt]\\[-8pt]
&&\eqntext{\mbox{where } \tilde Z_t(\omega):= \psi(t, \omega_{t_1\wedge
t},
\ldots, \omega_{t_n\wedge t}) \qquad\mbox{for all }(t,\omega)\in\L.}
\end{eqnarray}
Next, for some $\tilde h:= \tilde h_\e>0$ which will be specified
later, denote
%
\begin{equation}
\label{Ze} Z^\e_t:= \frac{1}{\tilde h}\int
_{t-\tilde h}^t \tilde Z_{\t\vee
s}\,ds \qquad\mbox{for } t
\in[\t, T].
\end{equation}
By choosing $\tilde h>0$ small enough (which may depend on $h_\e$), we have
%
\begin{equation}
\label{Zeest} \dbE^{\dbP_0} \biggl[\int_\t^T
\bigl|Z^\e_t - Z_t\bigr|^2\,dt \biggr] <2h.
\end{equation}

Now denote
\[
\d Z^\e:= Z^\e- Z,\qquad \d X^\e:=
X^\e- X.
\]
Then
\[
\d X^\e_t = \d X^\e_\t- \int
_\t^t \bigl[\a_s \d
X^\e_s +\bigl\langle\b, \d Z^\e\bigr
\rangle_s\bigr] \,ds + \int_\t^t \d
Z^\e_s \,dB_s,
\]
where $|\a|\le L_0$ and $\b\in\cU^{L_0}_t$. Denote
$
\G^\e_t:= \exp(\int_\t^t \a_s \,ds ).
$
We get
\[
\G^\e_t \d X^\e_t = \d
X^\e_\t- \int_\t^t
\G^\e_s\bigl\langle\b, \d Z^\e\bigr
\rangle_s \,ds + \int_\t^t
\G^\e_s \d Z^\e_s
\,dB_s.
\]
Then
\begin{eqnarray*}
0&\le& \sup_{\t\le t\le T} e^{-L_0t}\bigl|\d X^{\e}_t \bigr|
- e^{-L_0\t} \bigl|\d X^\e_\t\bigr| \le e^{-L_0\t}
\Bigl[\sup_{\t\le t\le T} \G^\e_t\bigl|\d X^{\e}_t
\bigr| - \bigl|\d X^\e_\t\bigr| \Bigr]
\\
&\le& \sup_{\t\le t\le T} \bigl|\G^\e_t\d
X^{\e}_t - \d X^\e_\t\bigr| =
\sup_{\t\le t\le T} \biggl|- \int_\t^t
\G^\e_s\bigl\langle\b, \d Z^\e\bigr
\rangle_s \,ds + \int_\t^t
\G^\e_s \d Z^\e_s \,dB_s
\biggr|
\\
&\le& C\int_\t^T\bigl|\d Z^\e_s\bigr|
\,ds + \sup_{\t\le t\le T} \biggl|\int_\t^t
\G^\e_s \d Z^\e_s \,dB_s
\biggr|.
\end{eqnarray*}
Thus
\begin{eqnarray*}
&&\dbP_0 \Bigl[\sup_{\t\le t\le T} e^{-L_0t}\bigl|X^{\e}_t
- X_t\bigr| \ge\e+e^{-L_0\t} |\tilde X_\t-
X_{\t}| \Bigr]
\\
&&\qquad = \dbP_0 \Bigl[\sup_{\t\le t\le T} e^{-L_0t}\bigl|X^{\e}_t
- X_t\bigr| - e^{-L_0\t} |\tilde X_\t- X_{\t}|
\ge\e\Bigr]
\\
&&\qquad\le \dbP_0 \biggl[C\int_\t^T\bigl|\d
Z^\e_s\bigr| \,ds + \sup_{\t\le t\le T} \biggl|\int
_\t^t\G^\e_s \d
Z^\e_s \,dB_s \biggr|\ge\e\biggr]
\\
&&\qquad\le \frac{C}{\e^2} \dbE^{\dbP_0} \biggl[ \biggl(\int
_\t^T\bigl|\d Z^\e_s\bigr| \,ds
\biggr)^2 + \sup_{\t\le t\le T} \biggl|\int_\t^t
\G^\e_s \d Z^\e_s \,dB_s
\biggr|^2 \biggr]
\\
&&\qquad\le \frac{C}{\e^2} \dbE^{\dbP_0} \biggl[\int_\t^T\bigl|
\d Z^\e_s\bigr|^2 \,ds \biggr] \le\frac{Ch}{\e^2},
\end{eqnarray*}
thanks to (\ref{Zeest}). Now set $h:= \frac{\e^3}{C}$, and we prove
(\ref{BBest1}).\vadjust{\goodbreak}

Finally, by (\ref{Ze}) and (\ref{tildeZ}) we have
\[
\th^\e_s= \frac{1}{\tilde h} [\tilde Z_s -
\tilde Z_{(s-\tilde
h)\vee\t}],\qquad s\in[\t, T].
\]
Fix $\omega\in\O$ and set $t:= \t(\omega)$. For each $\hat
\omega\in\hat\O^t$, set $\bar\omega:= \omega\otimes_t \hat
\omega\in\hat\O$, and define
\begin{eqnarray*}
\hat Z^{t,\omega}_s(\hat\omega)&:=& \psi(s, \bar
\omega_{t_1\wedge s},\ldots, \bar\omega_{t_n\wedge s}), \\
\bigl(\th^\e\bigr)^{t,\omega}_s(\hat\omega)&:=&
\frac{1}{\tilde h} \bigl[\hat Z^{t,\omega}_s(\hat\omega) - \hat
Z^{t,\omega}_{(s-\tilde h)\vee
t}(\hat\omega)\bigr],\qquad s\in[\t, T].
\end{eqnarray*}
Then we can easily check that $(\th^\e)^{t,\omega}$ satisfies (\ref
{Zproperty}).
\end{pf}
%
\begin{lem}
\label{lem-BB2}
Assume Assumption \ref{assum-fg} holds. Let $x\in\dbR$ and $Z$ be
$\dbF$-pro\-gressively measurable such that $\dbE^{\dbP_0}[\int_0^T
|Z_s|^2\,ds] <\infty$. For any $\e>0$, there exists $\dbF
$-progressively measurable process $Z^\e$ and an increasing sequence
of $\dbF$-stopping times $0=\t_0\le\t_1\le\cdots\le T$ such
that:

\begin{longlist}
\item
It holds that
%
\begin{equation}
\label{BBest2} \sup_{0\le t\le T} \bigl|X^{\e}_t -
X_t\bigr| \le\e, \qquad\dbP_0\mbox{-a.s.},
\end{equation}
where $X, X^\e$ are the solutions to the following ODEs with random
coefficients
%
\begin{eqnarray}
\label{Xe2} X_t &=& x - \int_0^t
f(s,B, X_s, Z_s) \,ds + \int_0^t
Z_s \,dB_s,
\nonumber
\\
X^\e_t &=& x - \int_0^t
f\bigl(s,B, X^\e_s, Z^\e_s\bigr)
\,ds \\
&&{}+ \int_0^t Z^\e_s
\,dB_s,\qquad 0\le t\le T,\dbP_0\mbox{-a.s.}\nonumber
\end{eqnarray}

\item For each $i$, $\th^\e_t:= \frac{d}{dt} Z^\e_t$ exists for
$t\in[\t_i, \t_{i+1})$, where $\th^\e_\t$ is understood as the
right derivative. Moreover, there exists $\tilde\th^{i,\e}$ on $[\t_i,
T]$ such that $\tilde\th^{i,\e}_t = \th^\e_t$ for $t\in[\t_i, \t
_{i+1})$, and for each $\omega\in\O$, $(\th^{i,\e})^{\t
_i(\omega), \omega}$ satisfies (\ref{Zproperty}) with \mbox{$t:= \t_i(\omega)$}.

\item For $\dbP_0$-a.s. $\omega\in\O$, for each $i$, $\t_i < \t_{i+1}$
whenever $\t_i < T$, and the set $\{i\dvtx  \t_i(\omega) < T\}$
is finite.
\end{longlist}
\end{lem}
\begin{pf} Let $\e>0$ be fixed, and set $\e_i:= 2^{-i-2}e^{-L_0T}\e
$, $i\ge0$. We construct $\t_i$ and $(Z^{i,\e}, X^{i,\e})$ by
induction as follows.

First, for $i=0$, set $\t_0:=0$. Apply Lemma \ref{lem-BB1} with
initial time $\t_0$, initial value $x$ and error level $\e_0$, we can
construct $Z^{0,\e}$ and $X^{0,\e}$ on $[\t_0, T]$ satisfying the
properties in Lemma \ref{lem-BB1}. In particular,
\[
\dbP_0 \Bigl[\sup_{\t_0\le t\le T} e^{-L_0t}\bigl|X^{0,\e}_t
- X_t\bigr| \ge\e_0 \Bigr] \le\e_0.
\]
Denote
%
\begin{equation}
\label{tau1} \t_1:= \inf\bigl\{ t\ge\t_0\dvtx
e^{-L_0t}\bigl|X^{0,\e}_t - X_t\bigr| \ge
\e_0 \bigr\} \wedge T.
\end{equation}
Since $X$ and $X^{0,\e}$ are continuous, we have $\t_1 > \t_0$,
$\dbP_0$-a.s. We now define
\[
Z^\e_t:= Z^{0,\e}_t,\qquad
X^\e_t:= X^{0,\e}_t,\qquad t\in[
\t_0, \t_1).
\]

Assume we have defined $\t_i$, $Z^\e, X^\e$ on $[\t_0, \t_i)$ and
$X^{i-1,\e}$ on $[\t_{i-1}, T]$. Apply Lemma \ref{lem-BB1} with
initial time $\t_i$, initial value $X^{i-1,\e}_{\t_i}$ and error
level $\e_i$, we can construct $Z^{i,\e}$ and $X^{i,\e}$ on $[\t_i,
T]$ satisfying the properties in Lemma \ref{lem-BB1}. In particular,
\[
\dbP_0 \Bigl[\sup_{\t_i\le t\le T} e^{-L_0t}\bigl|X^{i,\e}_t
- X_t\bigr| \ge\e_i + e^{-L_0\t_i} \bigl|X^{i-1,\e}_{\t_i}
- X_{\t_i}\bigr| \Bigr] \le\e_i.
\]
Denote
\[
\t_{i+1}:= \inf\bigl\{ t\ge\t_i\dvtx  e^{-L_0t}\bigl|X^{i,\e}_t
- X_t\bigr| \ge\e_i + e^{-L_0\t_i} \bigl|X^{i-1,\e}_{\t_i}
- X_{\t_i}\bigr| \bigr\} \wedge T.
\]
Since $X$ and $X^{i,\e}$ are continuous, we have $\t_{i+1} > \t_i$
whenever $\t_i<T$. We then define
\[
Z^\e_t:= Z^{i,\e}_t,\qquad
X^\e_t:= X^{i,\e}_t,\qquad t\in[
\t_i, \t_{i+1}).
\]

From our construction we have $\dbP_0(\t_{i+1} < T) \le\e_i$. Then
\[
\sum_{i=0}^\infty\dbP_0(
\t_{i+1} < T) \le\sum_{i=0}^\infty
\e_i <\infty.
\]
It follows from the Borel--Cantelli lemma that the set $\{i\dvtx  \t_i(\omega
) < T\}$ is finite, for $\dbP_0$-a.s. $\omega\in\O$,
which proves (iii).

We thus have defined $Z^\e, X^\e$ on $[0,T]$, and the statements in
(ii) follow directly from Lemma \ref{lem-BB1}. So it remains to prove
(i). For each $i$, by the definition of $\t_i$ we see that
\[
e^{-L_0\t_{i+1}}\bigl|X^{\e}_{\t_{i+1}} - X_{\t_{i+1}}\bigr| \le
\e_i + e^{-L_0\t_i}\bigl|X^{\e}_{\t_i} -
X_{\t_i}\bigr|,\qquad\dbP_0\mbox{-a.s.}
\]
Since $X^\e_{\t_0}= X_{\t_0}=x$, by induction we get
\[
\sup_i e^{-L_0\t_i}\bigl|X^{\e}_{\t_i} -
X_{\t_i}\bigr| \le\sum_{i=0}^\infty
\e_i \le\sum_{i=0}^\infty2^{-i-2}e^{-L_0T}
\e= \frac12 e^{-L_0T}\e,\qquad\dbP_0\mbox{-a.s.}
\]
Then for each $i$,
\begin{eqnarray*}
\sup_{\t_i\le t\le\t_{i+1}} \bigl|X^{\e}_t - X_t\bigr| &\le&
e^{L_0T} \bigl[\e_i + \bigl|X^{\e}_{\t_i} -
X_{\t_i}\bigr| \bigr] \\
&\le& e^{L_0T} \biggl[2^{-i-2}
e^{-L_0T} \e+ \frac12 e^{-L_0T}\e\biggr] \le\e, \qquad\dbP_0
\mbox{-a.s.},
\end{eqnarray*}
which implies (\ref{BBest2}).\vadjust{\goodbreak}
\end{pf}
\begin{pf*}{Proof of Theorem \ref{thm-peron}} Without loss of
generality, we shall only prove $\overline u_0 = u^0_0$. Recall that
$(Y^0, Z^0)$ is the solution to BSDE (\ref{BSDE}).
Set $Z:= Z^0$ and $x:=Y^0_0$ in Lemma \ref{lem-BB2}, we see that $X =
Y^0=u^0$ and thus $X$ satisfies the regularity in Proposition \ref{prop-u0}.

From the construction in Lemma \ref{lem-BB2} and then by Lemma \ref
{lem-BB1} we see that $\tilde\th^{0,\e}_t:= \frac{d}{dt} Z^{0,\e
}_t$ exists for all $t\in[0,T)$ and satisfies (\ref{Zproperty}). Then
by Lemma \ref{lem-Xsmooth} we see that $X^{0,\e}\in C^{1,2}_b(\L)$
and $\cL X^{0,\e} = 0$. This implies that, for the $\t_1$ defined in
(\ref{tau1}), $\t_1(\omega) > \t_0$ for all $\omega\in\O$ and,
by Example \ref{eg-tau}, $\t_1\in\cT_+$.

For $i=1,2,\ldots\,$, repeat the above arguments and by induction we can
show that, for each $i$ and each $\omega\in\O$, $\t^{\t_i(\omega
),\omega}_{i+1}\in\cT^{\t_i(\omega)}_+$. Moreover, by Lemma \ref
{lem-BB2}, $\{i\dvtx  \t_i <T\}$ is finite, $\dbP_0$-a.s.

We now let $u^\e$ be the solution to the following ODE:
\[
u^\e_t = X^\e_0 +
e^{L_0T}\e- \int_0^t f\bigl(s,B,
u^\e_s, Z^\e_s\bigr) \,ds + \int
_0^t Z^\e_s
\,dB_s.
\]
For $i=0,1,\ldots\,$, by the construction of $Z^\e$ in Lemma \ref
{lem-BB2} and following the arguments in Lemma \ref{lem-Xsmooth}, one
can easily show that
%
\begin{equation}
\label{Lue=0} u^\e\in\bar C^{1,2}_{\dbP_0}
\bigl([0,T]\bigr) \quad\mbox{and}\quad \cL u^\e= 0.
\end{equation}
Moreover, note that
\[
u^\e_t - X^\e_t =
e^{L_0T} \e- \int_0^t \a_s
\bigl[u^\e_s - X^\e_s\bigr] \,ds,
\]
where $|\a|\le L_0$. By standard arguments one has
\[
\sup_{0\le t\le T}\bigl|u^\e_t - X^\e_t\bigr|
\le e^{2L_0T}\e\quad \mbox{and}\quad u^\e_T-
X^\e_T \ge e^{-LT}\bigl[u^\e_0
- X^\e_0\bigr] = \e.
\]
Therefore, by (\ref{BBest2}) and noting that $u^0$ is bounded, $u^\e$
is bounded and
\[
u^\e_T(\omega) \ge X^\e_T(
\omega) + \e\ge X_T(\omega) = Y^0_T(\omega) =
g(\omega)\qquad\mbox{for } \dbP_0\mbox{-a.s. }\omega.
\]
This, together with (\ref{Lue=0}), implies that $u^\e\in\overline
\cD(0,0)$. Then, by the definition of $\overline u$,
\[
\overline u_0 \le u^\e_0 =
X^\e_0 + e^{L_0T}\e\le u^0_0
+ \e+ e^{L_0T}\e.
\]
Since $\e$ is arbitrary, we obtain $\overline u_0 \le u^0_0$.
\end{pf*}



\printaddresses

\end{document}